\documentclass[12pt]{article}
\usepackage[english]{babel}
\usepackage{latexsym, amssymb}
\usepackage{amsmath, color}
\usepackage{amsfonts}

\newtheorem{lemma}{Lemma}
\newtheorem{theorem}{Theorem}

\begin{document}
%\newenvironment{proof}{Proof}[section]
%\def\theequation{\arabic{section}.\arabic{equation}}

%-------------------------------------------------------------------
\def\a{\alpha}   \def\g{\gamma}  \def\b{\beta}  \def\th{\theta}
\def\l{\lambda}  \def\w{\omega}   \def\W{\Omega} \def\e{\varepsilon}
\def\d{\delta}   \def\f{\varphi}  \def\D{\Delta}    \def\r{\rho}
\def\s{\sigma}   \def\G{\Gamma}   \def\L{\Lambda}
\def\pd{\partial} \def\ex{\exists\,}  \def\all{\forall\,}
\def\dy{\dot y}   \def\dx{\dot x}   \def\t{\tau}
\def\bx{\bar x}   \def\bu{\bar u}   \def\bw{\bar w}  \def\bz{\bar z}
\def\p{\psi} \def\bh{\bar h}
\def\calO{{\cal O}}   \def\calP{{\cal P}}  \def\calR{{\cal R}}
\def\calN{{\cal N}}   \def\calL{{\cal L}}   \def\calQ{{\cal Q}}
\def\hx{\hat x}  \def\hw{\hat w}   \def\hu{\hat u}  \def\hv{\hat v}
\def\hz{\hat z}   % \def\ht{\hat t}

\def\empty{\O}     %\def\empty{\emptyset}
\def\lee{\;\le\;}   \def\gee{\;\ge\;}
\def\le{\leqslant}  \def\ge {\geqslant}
\def\Proof{{\bf Proof.}\quad}

\def\ov{\overline}  \def\ld{\ldots}
\def\tl{\tilde}   \def\wt{\widetilde}  \def\wh{\widehat}
\def\lan{\langle}  \def\ran{\rangle}

\def\hsp{\hspace*{20mm}} \def\vs{\vskip 8mm}
\def\bs{\bigskip}     \def\ms{\medskip}   \def\ssk{\smallskip}
\def\q{\quad}  \def\qq{\qquad}    \def\vad{\vadjust{\kern-5pt}}
\def\np{\newpage}   \def\noi{\noindent}
\def\pol{\frac12\,}    \def\const{\mbox{const}\,}
\def\R{\mathbb{R}}   \def\ctd{\hfill $\Box$}  \def\inter{\rm int\,}

\def\dis{\displaystyle}
\def\lra{\longrightarrow}  \def\iff{\Longleftrightarrow}
\def\beq{\begin{equation}\label}  \def\eeq{\end{equation}}
\def\bth{\begin{theorem}\label}   \def\eth{\end{theorem}}
\def\ble{\begin{lemma}\label}     \def\ele{\end{lemma}}
\def\begar{\begin{array}} \def\endar{\end{array}}

\newcommand{\vmin}{\mathop{\rm vraimin}}
\newcommand{\vmax}{\mathop{\rm vraimax}}

\newcommand{\blue}[1]{\textcolor{blue}{#1}}
\newcommand{\red}[1]{\textcolor{red}{#1}}
\newcommand{\gr}[1]{\textcolor{green}{#1}}

\def\bls{\baselineskip}   \def\nbls{\normalbaselineskip}
%\bls =1.1\nbls
%----------------------------------------------------------------

\begin{center}
{\Large\bf On the proof of Pontryagin's maximum principle\\[6pt]
by means of needle variations\footnote{This paper is submitted to
an issue dedicated to 80th birthday of professor V.M. Tikhomirov.}
%Vladimir Mikhailovich Tikhomirov.}
}
\\[8mm]
{\large A.V. Dmitruk, \,\, N.P. Osmolovskii} \\[10pt]
%{\large А.В. Дмитрук\footnotemark[1], \quad Н.П. Осмоловский \footnotemark[2]}
%\\[20pt]
%{\large Moscow State University}\\[4pt]
%\footnotemark[1]{avdmi@cemi.rssi.ru}\,\,
%\footnotemark[2]{osmolovsky@uph.edu.pl}\\[20pt]
\end{center}
\bs

\small
{\bf Abstract}.  We propose a proof of the maximum principle for the
general Pontryagin type optimal control problem, based on packages of needle
variations. The optimal control problem is first reduced to a family of
smooth finite-dimensional problems, the arguments of which are the widths
of the needles in each packet, then, for each of these problems, the standard
Lagrange multipliers rule is applied, and finally, the obtained family
of necessary conditions is "compressed" in one universal optimality
condition by using the concept of centered family of compacta.
\ms

{\bf Keywords:}\, Pontryagin maximum principle, package of needle variations,
operator of terminal state value, smooth finite dimensional problem,
Lagrange multipliers, finite-valued maximum principle, centered family
of compacta.
\bs
\normalsize

\section{Introduction}

As is known, the original proof of the maximum principle (MP) as a necessary
condition for a strong minimum, given in the book of Pontryagin et al \cite{Pont},
was very complex and did not fit into the framework of classical methods
of the calculus of variations. Later, efforts of many mathematicians were aimed
at rethinking this proof and extending the MP to new classes of problems.
This was done in a large number of works, both Russian and international;\,
see, e.g., \cite{DM65} -- \cite{MDO} (of course, this list in no way pretends
to be complete).\, On the other hand, the question of the simplest and most
transparent proof of the MP for the {\em classical problem of Pontryagin type\,}
was also a subject of attention to specialists. This is especially important
from the educational and methodical point of view in order of teaching the
optimal control in both mathematical and engineering faculties of universities.
\ssk

The basic idea in the proof of MP is to introduce a certain class (family)
of control variations along which {\em it is possible to differentiate.\,}
For control systems of ODEs, the following classes of variations are known:\,
a) uniformly small variations,\, b) needle variations and their "packages",\,
c) the so-called $v- $change of time, and\, d) sliding mode regimes.
\ssk

The uniformly small variations correspond to the {\it weak}\, minimum and lead
to the {\em stati\-ona\-rity condition,} i.e. to the Euler-Lagrange equation.
It is a very important result that can serve as a starting point for further
investigations of the given process for optimality, but it is far from the MP.
For obtaining conditions of the strong minimum, these variations are
obviously not enough.
\ssk

The simplest example of\, "non-small" variation is the {\it needle variation}\,
of the control, consisting in the change of the optimal control by an
arbitrary admissible value of the control only on a small time interval.
The width of this interval is a parameter of variation, with respect to
which one can differentiate. If the right endpoint of trajectory is free
of any restrictions (the so-called {\it free endpoint problem}), then
the cost functional on this one-parametric family of trajectories must
have a minimum at zero, hence its derivative w.r.t. the needle's width
must be nonnegative. This easily yields the conditions of MP.\,
However, if the problem includes constraints on the right endpoint of
trajectory, such a variation may generate a trajectory that would not,
in general, satisfy them. In order to guarantee that the endpoint of
varied trajectories hit the given constraints, the family of variations
should be rich enough, and therefore one has to consider not just one,
but a finite number of needle variations together, the so-called {\it
packet}\, of needle variations, whose parameter is the collection of
widths of the needles, independent of each other.

Needle variations and their packages have a "genetic flaw" consisting
in that the widths of the needles are by definition nonnegative, so the
functions obtained with their help are  defined only on the nonnegative
orthant of a finite-dimensional space (more precisely, on its intersection
with a neighborhood of the origin), and the further study of the obtained
accessory problem corresponding to the given package of needles requires
to perform a number of specific constructions that go beyond the standard
courses of analysis and differential equations. This immediately heavies
the proof of MP, making it almost unacceptable for most of listeners.
Therefore, the proofs in the courses of optimal control are often limited
to the case of a free endpoint problem, where, as is known, it suffices
to use only one needle and no additional constructions are needed.
To present the proof in problems {\em with constraints on the endpoints
of trajectory}, it is desirable to have the accessory problem in the form
of a smooth problem of more or less standard (previously studied) type.
\ssk

One of the techniques that leads to a smooth  problem is the so-called
{\em $v$-change of time}, proposed by Dubovitskii and Milyutin in \cite{DM65}
(see also \cite{Gir, IT, MDO, OPU, develop}), consisting in the introduction
of a new time $\t$ and interpretation of the original time $t=t(\t)$ as a
state variable satisfying the equation
$$
\frac {dt}{d\t}\;=\; v(\t), \qq \mbox{where}\q v(\t)\ge 0
$$
is one more  control. Note that this change of time can be not one-to-one
(and this is important!).\, The  simplest version of this change, where
the dependence of the "old"\, time on the "new" one is piecewise linear
(i.e., $v(\t)$ is piecewise constant), was described in detail in \cite{MDO, OPU}.\,
The nontriviality of such a change is that the small variations of the
control $v$ lead to needle variations (or close to them) of the original
control $u,$ so that in fact $v-$change of time represents such a parametrization
of needle variations that admits a natural smooth extension to the "negative
widths". This trick is very nice, however, its usage requires, though
very transparent in ideas, technically rather cumbersome and not completely
standard constructions, and the teaching experience shows that such a
proof of MP\,  still is perceived by students with difficulty.
\ssk

Another way to prove the MP, also proposed by Dubovitskii and Milyutin,
is based on variations of {\em sliding modes type}. (The sliding modes
themselves were introduced by Gamkrelidze in \cite{Gam} to prove the existence
of solution to the problem.)\, This method, also elegant, allows far-reaching
generalizations (e.g., to problems with mixed constraints, see
\cite{VNIISI, MDO, develop}), but requires the usage of rather nontrivial
(and nonstandard) facts of linear and nonlinear functional analysis, which
makes it hardly justified in application to classical Pontryagin problem,
especially in a regular lecture course. \ms

Let us get back to the needle variations as the simplest class of "non-small"
variations.  Over the years, in the seminar at the Department of General
Control Problems\, of Mechanics and Mathematics Faculty of Moscow State
University, led by V.M. Tikho\-mi\-rov, it was intensively discussed the
program of the lecture course "Calculus of Variations and Optimal Control"
(CVOC), including various schemes of the proof of MP.\,  The main focus
was on the most simple description of the cone of state endpoint variations
generated by the packages of needle variations of the control.
To our opinion, here the most successful construction was proposed by
Magaril-Il'yaev \cite{Mag}.\,
The proof of MP itself follows in \cite{Mag} the scheme of the original
proof from the book \cite{Pont} (and also \cite{Bolt}). If this cone is
the full space, it is also the same for some fixed packet of needles, and
then the corresponding nonlinear operator defined on the nonnegative
orthant of needles widths in this packet possesses the so-called property
of {\it covering at the point} (the image of any neighborhood of the origin
contains a neighborhood of the optimal point), which leads to a contradiction
with the optimality of the reference trajectory.
If, otherwise, this cone is not the full space, the application of the
finite-dimensional separation theorem easily gives the MP. \,
The proof of the covering on a cone relies in \cite{Mag} either on a cone-type
modification of {\it the Newton method} (in the case of piecewise continuous
control), or on a rather fine and nontrivial topological fact, {\it the
Brower fixpoint theorem} (in the case of measurable control). None of these
facts is presented in standard courses, especially for engineers.
(In \cite{Pont, Bolt}, this proof relied on even more specific topological
notion of {\it intersection index}\, or {\it the Sperner lemma}, which
are not lectured even for mathematicians.)\, Moreover, only the problem
with fixed endpoints is considered in \cite{Mag}.
\ssk

Another way to prove the MP by needle variations, which seems to us more
preferable, and to which the present paper is devoted, is the passage to
a standard finite-dimensional problem, the arguments of which would be the
widths of the needles in a given package. This idea is not new (see eg.
\cite{ATF, Mich, Car} and a recent paper \cite{Kor}), but its implementation
requires to determine an extension of the needles to the negative widths.
(As was already said, an alternative to such extension is the above-mentioned
$v-$change of time.)\, In \cite[Sec. 4.2]{ATF}, it was proposed an explicit
extension of the operator of the terminal state value (corresponding to
a given package of needles) to a full neighborhood of the origin in the
space of needles' widths, but the construction of this extension is
realized rather cumbersomely, and therefore it still hardly can be successfully presented in a regular course.
(In papers \cite{Mich, Car}, the required extension was not discussed,
so, the proposed proofs of MP were, in fact, not complete.)
\ssk

Somewhen in the mid 2000-s, at the above mentioned seminar of V.M. Tikhomirov,
the authors of the present paper proposed to use a possibility of extension
of {\em an arbitrary}\, smooth function, defined on the intersection of
nonnegative orthant $\R^k_+$ with a neighborhood of zero, to a full neighborhood
of zero, preserving its smoothness. Of course, one can use here the classical
theorem of Whitney on the extension of a smooth function from an arbitrary
compact set to its neighborhood (see \cite{Malg}), but again, it would hardly
be reasonable to apply, in this simple case, such a difficult theorem, even
formulation of which is rather complicated and requires additional notions.
We proposed a simple way of such extension for an arbitrary function
defined on $\R^k_+$ and {\it strictly differentiable}\, only at zero.
This is quite enough for our purposes.
(If a function is {\it continuously differentiable}\, everywhere in a
neighborhood of zero in $\R^k_+\,,$ which is also enough for us, then
one can extend it to a full neighborhood of zero by the coordinate-wise
method, that is reduced to an obvious one-dimensional case. This method
was used in the paper \cite{Kor}, which is the first known to us publication
with an explicit construction of the required extension;\, see details
below.)\,

With this extension, no matter how it is constructed, there is no need to
define the needles of negative widths, nor to introduce and describe the
cone of endpoint variations, and the proof of MP is as follows.  For any
packet of needle variations we obtain an auxiliary (or "accessory") {\em
smooth finite-dimensional problem}\, in the space of needles' widths of
the given packet, in which the origin is a point of local minimum.
Applying the standard Lagrange multipliers rule, we obtain a "finite-valued"
MP, in which the maximality of the Pontryagin function holds {\em for
the given set of needles}. Thus, we have {\em a family of finite-dimensional
problems}\, corresponding to all possible packets of needles, and in each
problem we obtain its finite-valued MP.\, The tuples of Lagrange multipliers
in each of these problems is a finite-dimensional compact set (a simple
fact), and the family of these compacta turns out to be {\it centered}
(i.e., any finite number of them has a nonempty intersection). Therefore,
all this family also has a nonempty intersection, and any element of the
latter (i.e. a tuple of Lagrange multipliers) guarantees the "universal"
MP, in which the maximality condition holds {\it for any set of needles},
i.e. for any time instants and any admissible values of control.
\ssk

This proof relies only upon standard facts of mathematical analysis for
the first year students, differential equations, and the Lagrange multipliers
rule for the smooth finite-dimensional problem with equality and inequality
constraints (which is assumed to be given in the course of CVOC before the
presentation of MP).\, Like in other proofs, we also use a topological
fact, the finite intersection property of a family of compacta, but this
is one of equivalent definitions of compactness, which hopefully should
not present any difficulties in teaching, especially in the
finite-dimensional case.  \ssk

Actually, our proof follows the same scheme as in the book \cite{ATF};\,
it differs only in a more simple way of extension of functions on a packet
of needles and by a rather more simple presentation of differential properties
of these functions.  From that in \cite{Mag} it differs by that, instead
of description of the terminal cone and proving a theorem on the cone-type
covering property, we pass to a finite-dimensional problem, apply the
Lagrange multipliers rule, and then use the property of a centered
family of compacta. In our opinion, this way of proving MP for the classical
Pontryagin type problem (in its most general setting) is more simple in
some aspects than those commonly used in the lecture courses, and so it
deserves attention. \ssk

Note also that the scheme of proving MP based on the introduction of
a family of so-called {\it accessory}\, smooth problems, on writing out
the stationarity conditions in each of them, and then on using a centered
family of compacta (or, more general, of a projective family of compacta),
was proposed by A.Ya. Dubovitskii and A.A. Milyutin (see \cite{DM69, Mil70})
and was effectively applied not only to the Pontryagin type problem but
also to a more general class of problems that include both state and mixed
state-control constraints.  In this paper we do not consider these
generalizations, referring the reader to the works \cite{VNIISI, MDO, develop}.
\ms

Let us pass to the detailed exposition.\\

%%-------------------------
\section{General  problem of Pontryagin type, formulation of the
maximum pronciple}

On a time interval  $[t_0,t_1],$  not fixed a priori, we  consider the
following optimal control problem, which  will be  called {\it the general}
(or {\it canonical}) {\it Pontryagin type problem}:
\begin{equation}\label{4}
\dot x(t)=f(t,x(t),u(t)), \qquad u(t)\in U,
\end{equation}
\begin{equation}\label{3}
K(t_0,x(t_0),t_1,x(t_1))=0,
\end{equation}
\begin{equation}\label{2}
F(t_0,x(t_0),t_1,x(t_1))\le 0,
\end{equation}
\begin{equation}\label{1}
J =\; F_0(t_0,x(t_0),t_1,x(t_1))\to \min.
\end{equation}
Here $x\in\R^n$ is the state variable, $u\in\R^r$ is the control;\,
the endpoint functions  $K,\, F$ have dimensions  $d(K),\, d(F),$ respectively.
(Such notation for the dimensions of vectors is convenient since it allows to
"save" letters and does not require to remember the dimensions.)

We assume that the functions $F_0\,,F,\, K$ are of class $C^1,$ and $f$ is
continuous together with its derivatives  $f_t$ and $f_x\,.$ The set
$U \subset \R^r$ is arbitrary. For brevity, problem (\ref{4})--({\ref{1})
is called {\em problem} $A$.\ssk

{\bf Remark 1.}\, To be more precise, the properties of function $f$ can
be assumed to hold on the set ${\cal Q}\times U,$ where ${\cal Q}$ is an
open set in $\R^{1+n},$ and the properties of functions $F_0\,, F,\, K$
to hold on an open set ${\cal P}\subset\R^{2n+2}.$  As a rule, it is always
assumed without explicit indication of these sets. We also will not be
distracted to these minor details.  \ssk

{\bf Remark 2.}\, As was mentioned in the book \cite{Pont}, one need not
necessarily assume the control set $U$ to lie in a finite-dimensional space;\,
in general, it can be any Hausdorff topological space. All the below arguments
remain then unchanged. However, this generalization would distract the reader
(and listener) to inessential issues, while it hardly is necessary for
applications. Therefore, like in \cite{Pont}, we assume that $U \subset \R^r.$
\ms

The solution of problem A\, is sought in the class of absolutely continuous
functions  $x(t)$ and measurable bounded functions $u(t).$ A pair of functions
$w(t)=(x(t), u(t))$ together with a segment $[t_0, t_1]\,$ of their definition
is called\, {\it a process}\, of the problem.\\
A process is called {\em admissible}, if it satisfies all the constraints
of the problem. Conditions (\ref{4}) are assumed to hold almost everywhere.\,
As usual, we say that an admissible process
$\hw(t)=(\hx(t), \hu(t))\mid\, t\in [\hat t_0, \hat t_1],$ provides the
{\it strong minimum}\, if  there exists an $\e>0$ such that $J(w)\ge J(\hw)$
for all admissible processes $w(t)=(x(t),u(t))\mid\, t\in [t_0, t_1],$
satisfying the following conditions: \vad
$$
|t_0-\hat t_0|<\varepsilon, \qquad |t_1-\hat t_1|<\varepsilon,\quad
$$ $$
|x(t)-\hat x(t)|<\varepsilon  \qquad  \all t\in [t_0,t_1]\cap [\hat t_0,\hat t_1].
$$ \ssk

We now state\, {\bf the Pontryagin maximum principle},\, a necessary condition
for a strong minimum in problem A.\, Introduce {\em the Pontryagin
function}  \vad
$$
H(\p_x,t,x,u)\;=\; \psi_x f(t, x, u),
$$
where  $\psi_x$ is a row vector of dimension $n$ (the dependence of $H$ on
$\psi_x$ will be sometimes omitted), and {\em the endpoint Lagrange function}
$$
l(t_0,x_0,t_1,x_1)=\; (\alpha_0 F_0 + \alpha F + \beta K)(t_0,x_0,t_1,x_1),
$$
where $\alpha_0$ is a number, and $\alpha,\, \beta$ are row vectors of
the same dimensions as $F,\,K,$ respectively (we omit the dependence of
$l$ on $\alpha_0\,,\, \alpha,\,\beta$). \ms

Let  $w=(x(t), u(t))\mid\,  t\in [t_0, t_1],$ be an admissible process in
problem A.\, We say that it satisfies {\em the Pontryagin maximum principle}\,
if there exist a number  $\alpha_0\,,$ row vectors  $\alpha\in\R^{d(F)},$
$\beta\in \R^{d(K)},$ and  absolutely continuous functions  $\psi_x(t),\; \psi_t(t)$
of dimensions  $n,\,1,$ respectively (where  $x$ and $t$ are the subscripts,
not the notation of derivatives), such that
\begin{itemize} \item[(i)] $\alpha_0\ge0$, $\alpha\ge0$;

\item[(ii)] $\alpha_0+|\alpha|+|\beta|>0$;

\item[(iii)] $\alpha F(t_0,x(t_0),t_1,x(t_1))=0$;

\item[(iv)] $-\dot\psi_x(t)=H_x(\psi_x(t),t,x(t), u(t)),\qq
-\dot \psi_t(t)=H_t(\psi_x(t),t,x(t), u(t)),$

\item[(v)]

$\psi_x(t_0)= l_{x_0}(t_0,x(t_0),t_1,x(t_1)),\qquad
\psi_x(t_1)= -l_{x_1}(t_0,x(t_0),t_1,x(t_1)); \\[6pt]
\psi_t(t_0)= l_{t_0}(t_0,x(t_0),t_1,x(t_1)),\qquad
\psi_t(t_1)= -l_{t_1}(t_0,x(t_0),t_1,x(t_1)),$

\item[(vi)] $H(\psi_x(t),t,x(t), u(t)) + \psi_t(t) =0$ \q\,
for a.a. $t\in[t_0,t_1],$

\item[(vii)] $H(\psi_x(t),t,x(t), u')+ \psi_t(t) \le 0$\q
for all\, $t\in[t_0,t_1]$ and all\, $u'\in U$.

\end{itemize}

Conditions $(i)-(v)$ are called the nonnegativity, nontriviality, complementary
slackness, adjoint equations, and transversality conditions, respectively.\,
Con\-di\-tion $(vi)$ has not yet a standard name;\, "in working order"
we call it {\em the energy evolution law}, since $(vi)$ and the adjoint
equation for $\p_t$ yields the equation for the function $H,$ which in
mechanical problems is usually regarded as the energy of the system:
$$
\dot H \;=\; H_t\, \q\; \mbox{or} \q\; \frac{dH}{dt}\;=\; \frac{\pd H}{\pd t}\;.
$$
(If the control system is time-independent, i.e. $f = f(x,u),$ we get the
{\em energy conservation law}: $\;\dot H =0,$ i.e. $H = \const.$) \ssk

Conditions  $(vi)$ and $(vii)$ imply {\it the maximality condition}\, for
the Pontryagin function:  \vad
$$
\max_{u' \in U}\, H(\psi_x(t),t,x(t), u') \; =\; H(\psi_x(t),t,x(t), u(t))
\q\, \mbox{for almost all}\;\; t\in[t_0,t_1],
$$
which gave to the whole set of conditions $(i)-(vii)$ the name of
{\em Pontryagin maximum principle.}
\ms

Note that notation $\psi_x(t)$ and $\psi_t(t)$ for the adjoint variables
was proposed by A.Ya. Dubovitskii and A.A. Milyutin. The convenience of
such notation is quickly clarified in solving concrete problems with
multiple state variables.
\ms

{\bf Remark 3.}\, One can show that the equation for the function $\p_t$
(i.e. for $-H)$ follows from the other conditions of MP, i.e. it is not
independent. (Below, we will show it for the time-independent problem.)\,
Nevertheless, including it in the set of conditions of MP\, is justified
not only by its relation with the principles of mechanics, but also by
the fact that, in many problems, this condition is very convenient to use
directly, in its "ready form", not deriving it from other conditions of MP.
\ms

Necessary conditions for a strong minimum are given by the following theorem.

\begin{theorem}\label{th1}
If a process  $\hw=(\hx(t), \hu(t))\mid  t\in [\hat t_0, \hat t_1])$
provides the strong minimum in problem $A,$ then it satisfies the Pontryagin
maximum principle.
\end{theorem}

First we prove this theorem for the case when the  time interval is fixed
and the control system is autonomous. (The latter, however, is not necessary.)
\bs

%%-------------------------------
\section{Maximum principle for the problem on a fixed time interval}

Consider the following problem $B$:
\begin{equation}\label{B4}
\dot x(t)= f(x(t),u(t)), \qquad u(t)\in U, \qq t\in[t_0,t_1],
\end{equation}
\begin{equation}\label{B3}
K(x(t_0),\,x(t_1))=0, \q
%\end{equation}
%\begin{equation}\label{B2}
F(x(t_0),\,x(t_1)) \le 0, \q
\end{equation}
\begin{equation}\label{B1}
J =\; F_0(x(t_0),\,x(t_1))\to \min.
\end{equation}
Here the interval  $[t_0,t_1]$ is fixed, the functions  $F_0\,,\, F$ and $K$
are continuously differentiable, and the function  $f$ is continuous together
with its derivative $f_x$.\, The minimum is sought among all pairs (processes)
$w=(x,u) \in W: = AC^n[t_0,t_1]\times L_\infty^r[t_0,t_1],$ where $AC^n[t_0,t_1]$
is the space of absolutely continuous functions of di\-men\-sion $n,$
and $L_\infty^r[t_0,t_1]$ is the space of measurable bounded functions
of di\-men\-sion $r$.
\ssk

By definition, an admissible  process  $\hw=(\hx,\hu)$ provides {\em the
strong minimum\,} if there exists $\e >0$ such that $J(w)\ge J(\hat{w})$
for all admissible pairs  $w=(x, u)$ satisfying the condition
$||x-\hat x||_C <\e.$ In other words,  the strong minimum is a local
minimum in the space $W$ with respect to the seminorm $||w||' =||\,x||_C\,.$
\ms

As in above, we introduce the Pontryagin function $H(\p,x,u)\,=\, \psi f(x, u),$
where  $\psi$ is an $n-$dimensional row vector (in the general problem A,
it was denoted by $\psi_x),$ and the endpoint Lagrange function
$l(x_0,\,x_1)=\; (\alpha_0 F_0 + \alpha F + \beta K)(x_0,\,x_1).$
\ssk

For an admissible process $w=(x, u)$ in problem B, the conditions of MP\,
consist in the following:\, there exist a number  $\alpha_0\,,$ row vectors
$\alpha\in\R^{d(F)},$ $\beta\in \R^{d(K)},$ an absolutely continuous function
$\psi(t)$ of dimension $n,$ and a constant $c\in\R,$ such that \vad

\begin{itemize} \item[(i')] $\alpha_0\ge 0$, $\alpha\ge 0$;

\item[(ii')] $\alpha_0+|\alpha|+|\beta|>0$;

\item[(iii')] $\alpha F(x(t_0),\,x(t_1))=0$;

\item[(iv')] $-\dot\psi(t)=\; H_x(\p(t), x(t), u(t)),$

\item[(v')]

$\psi(t_0)= l_{x_0}(x(t_0),\,x(t_1)),\qquad
\psi(t_1)= -l_{x_1}(x(t_0),\,x(t_1)),$

\item[(vi')] $H(\p(t), x(t), u(t)) =\; c$ \q\, almost everywhere on $[t_0,t_1],$

\item[(vii')] $H(\p(t), x(t), u')\lee c$ \q for all\,
$t\in[t_0,t_1]\;$ and all\, $u'\in U.$

\end{itemize}

Note that here, instead of the adjoint equation for $\psi_t\,,$ we write
the condition that the Pontryagin function is constant along the optimal
process (the ``energy conservation law'',\, as is expected in the
conservative system).
\ms

The following theorem holds.

\begin{theorem}\label{th2}\,\,
If a pair  $\hat w=(\hat x, \hat u)$ provides the strong minimum in
problem B, then it satisfies\, the maximum principle $(i')-(vii').$
\end{theorem}

The proof will be given here for the case where the control $\hat u(t)$
is piecewise continuous.  The general case, where the optimal control is
measurable, is technically more complicated and, as a rule, not
considered in the regular courses.
\ssk

First, we study some differential properties of the control system
\beq{xua1}
\dx = f(x,u)\qq % x(t_0)= x_0\,
\eeq
along its arbitrary solution  $(x(t),\,u(t))$ on the interval $[t_0,t_1].$
We will consider the controls uniformly bounded by some constant:
$|\,u(t)| \le M.$ From the ODE theory it is well-known (though it is also
easy to prove directly) that, if in addition $||u-\hu||_1$ is small and
the initial value $x_0=x(t_0)$ is close to  $\hx(t_0),$ then the solution
$x(t)$ of equation (\ref{xua1}) always exists on the entire interval
$[t_0,t_1]$ and is close to $\hx(t)$ in the norm of $C[t_0,t_1].$
\ssk

For a given process $(x(t),\,u(t))$ with a piecewise continuous control,
introduce
\ms

%%--------------------------------------
{\bf The concept of an elementary needle variation}.\,
Let us assume that the control $u(t)$ is left continuous at all points in
the semi-interval $(t_0,t_1]$ and right continuous at $t= t_0\,.$
Fix any point $\theta\in (t_0,t_1)$ and any value $v\in U.$
For each sufficiently small $\varepsilon\ge 0,$ define, on the interval
$[t_0,t_1],$ the new control
$$
u_\e(t)\,= \;
\left\{\begin{array}{l}
\; v,\quad\;\; \mbox{if}\quad t\in (\theta -\e,\,\theta),\\[4pt]
 u(t),\;\;\; \mbox{if}\quad t\notin (\theta -\e,\,\theta).
\end{array}\right.
$$
(If $\e =0,$ the semi-interval $(t_0,t_1]$ is empty and the control does
not change.) \\ Note that this control is still left continuous.
The family of control functions $u_\e(t)$ with $\e \to 0+$ is called
{\em an elementary needle variation}\, of the original control $u(t)$
with the parameters $(\theta, v).$ \ms

Let a vector $a\in\R^n$ be close to $ x_0= x(t_0).\,$ On the interval
$[t_0,t_1],$ define the state variable $x_\e(t)$ as the solution
to Cauchy problem for the control system (\ref{xua1}):
\beq{xea1}
\dx_\e =\, f(x_\e,\, u_\e), \qq x_\e(t_0) =a.
\eeq
If $\e\ge 0$ is small enough and   $a$ is close enough to  {$ x_0\,,$}
then, as was already said, the solution of problem (\ref{xea1}) exists
on the whole interval $[t_0,t_1],$ is unique, and continuously depends
on the pair $(a,\e)$ in the norm of $C[t_0,t_1].$ We accept this as an
established fact. The value of state variable  $x_\e(t)$ at the point
$t_1$ will be denoted by $P(a,\e).$
\ssk

Thus, we obtain a mapping
$$
P:\; (a,\e)\; \longmapsto\; x_\e(t_1) \in \R^n,
$$
which is defined and continuous on $\W(x_0) \times [0,\e_0),$ where
$\W(x_0)$ is a neighborhood of the point  $x_0 =x(t_0) \in \R^n\,$ and
$\e_0 >0.$ Obviously, $P( x_0,0)=  x(t_1)=: x_1$.
\ms

Let us show that, at the point $(a,\e)= (x_0,\,0),$ the mapping $P$ has
the derivative with respect to $a$ and the right derivative with respect
to $\e,$ which in the sequel will be denoted by $P_a(x_0,0)$ and
$P_\e^+(x_0,0),$ respectively. \ssk

As is known, to find the derivative with respect to $a$ in the direction
$\bar a,$ one has to solve {\it the equation in variations}\,
\beq{urvar1}
\dot{\bx} =\, f_x(x,\, u)\,\bx, \qq \bx(t_0)= \bar a.
\eeq
Then $P_a(x_0,0)\,\bar a\, =\, \bx(t_1).$  Note that if $x(t)$ slightly
varies with respect to the norm of the space $C[t_0,t_1],$ and $u(t)$ slightly
varies with respect to the norm of the space $L_1[t_0,t_1]$ preserving
the uniform boundedness with an a priori constant, then the matrix $f_x(x, u)$
not much varies with respect to the norm of $L_1[t_0,t_1],$ and then,
as is well known, the solution to equation (\ref{urvar1}) not much varies
in the norm of $C[t_0,t_1].$
\ms

To find the right derivative $P_{\e}^+(x_0,0)$ for a fixed $x_0\,,$
we will reason as follows. \linebreak
To each $\e> 0$, there corresponds a solution  $x_\e(t)$ to equation
(\ref{xea1});\,\, for $\e=0$ we have the unperturbed solution  $x(t).$
If $t<\theta,$ then, for sufficiently small $\e,$ both solutions are the
same at this point. Let us estimate the difference $\D x_\e = x_\e - x$
at the point $\theta.$ We  have
$$
\begin{array}{l}
x(\theta)\, =\; x(\theta-\e)\, +\, \e\,f(x(\theta), u(\theta)) + o(\e),  \\[6pt]
x_\e(\theta) =\; x(\theta-\e)\, +\, \e\,f(x(\theta), v)\; +\; o(\e),
\end{array}
$$ \\[2pt]
whence $\D x_\e(\theta) = \e\cdot \D f(\theta,v) + o(\e),\;$ where
$\D f(\theta,v) = f(x(\theta), v) - f(x(\theta), u(\theta)).$
\ms

Set  $\dis \bx(t) = \lim_{\e\to 0+}\, \frac{\D x_\e(t)}{\e}\,.\;\,$
Then  $x_\e(t) = x(t) + \e\,\bx(t) + o(\e).$ \\[8pt]
Moreover, to the left of $\theta,$ we obviously have $\bx(t) =0;\,$ at this
very point $\bx(\theta) = \D f(\theta,v),$ and on the interval $[\theta,\,t_1]$
both functions $x$ and $x_\e$ satisfy (\ref{xua1}) with the same unperturbed
control, but with different initial values at the point $\theta.$
Therefore, like before, $\bx$ satisfies on the interval $[\theta,\,t_1]$
the equation in variations
\beq{urbx1}
\dot{\bx} =\, f_x(x,\, u)\,\bx \qq
\mbox{with the initial value }\q \bx(\theta)= \D f(\theta,v),
\eeq
and then $\, P_{\e}^+(x_0,0) = \bx(t_1).$ \ms

Note that $P_\e^+(x_0,\,0)$ depends not only on the needle's parameters
$(\theta,v),$ but also on the base pair $(x(t), u(t))$ satisfying equation
(\ref{xua1}).\, If the base pair not much varies in the norm
$||x||_C + ||u||_1$ (with the uniformly bounded control), and the pair
$(\theta,\,u(\theta))$ also not much varies, then the matrix $f_x(x(t), u(t))$
will not much vary in the norm of $L_1\,,$ the initial time $\theta$
and the initial value $\bx(\theta)$ will not also much vary, hence the
solution to Cauchy problem (\ref{urbx1}) will have just a small change
in the norm of $C,$ and so, its terminal value $x(t_1)$ will have just
a small change.
\ms

%%-----------------------------------------------------
Thus, we proved the following

\ble{lem1}
At the point $(a,\e)= (x_0,\,0),$ the mapping $P$ has a derivative with
respect to $a$ and the right derivative with respect to $\e,$ which are
expressed by the formulas:
$$
\all \bar a \in \R^n \qq P_a(x_0,\,0)\,\bar a\, =\, \bx(t_1),
$$
where $\bx(t)$ is the solution of equation (\ref{urvar1}) on the interval
$[t_0,\,t_1],$ while the right derivative $\, P_{\e}(x_0,0) = \bx(t_1),$
where  $\bx(t)$ is the solution to Cauchy problem (\ref{urbx1}) on the
interval $[\theta,\,t_1].\;$
\ssk

Both these derivatives depend continuously on the pair $(x, u)$ in the
norm of the space $C \times L_1$ (under the uniform boundedness of
the control) and on the pair $(\theta,\,u(\theta)) \in \R^{1+r}.$
\ele
\ssk

Now, if the base control is taken to be $u_\e(t)$ for a small $\e>0,$
the initial state value is taken equal to a vector $a$ close to $x(t_0),$
and the needle variation is set at the point $\theta -\e,$ then the difference
$u_\e -u$ is small in the norm $L_1[t_0, t_1],$ the difference $x_\e -x$
is small in the norm $C[t_0, t_1],$ and in view of left continuity of the
function $u(t),$ the new base value $u_\e(\theta -\e) = u(\theta -\e)$
is close to the old one $u(\theta).$ Then the solution to Cauchy problem
(\ref{urbx1}) on the interval $[\theta -\e,\,t_1]$ restricted to the
interval $[\theta,\,t_1]$ is uniformly close to the old solution.
By these arguments, Lemma \ref{lem1} implies  \ms

{\bf Corollary.}\q {\em The operator $P$ corresponding to the elementary
needle variation at the point $\theta,$ for all initial values $a$ close
enough to $x(t_0)$ and all small enough $\e \ge 0,$ has at the point
$(a,\e)$ the derivative $P_a(a,\e)$ and right derivative $P_\e^+(a,\e)$
which continuously depend on the pair $(a,\e).$}
\ssk

The exact expressions of these derivatives will be needed for us only
at the point $(a,\e) = (x_0,\, 0);\,$ they are given above. \bs

\small
{\bf Remark 4.}\, In the case when the base control is measurable, this
Corollary is no longer true even if $\theta$ is {\it a Lebesgue point}\,
of the function $u(t).$ The matter is that, in general, the initial value
for equation (\ref{urbx1}) is given by the formula
$$
\bx(\theta) =\; f(x(\theta), v) -\; \lim_{\d \to 0}\,
\frac 1\d \int_{\theta -\d}^\theta f(x(\theta), u(\t))\,d\t,
$$
if this limit exists. At the very point $\theta$ it does exist and equals
$f(x(\theta), u(\theta)),$ so still $\bx(\theta)= \D f(\theta,v),$ but
at the shifted point $\theta -\e\,$ for $\e>0$ this limit can be essentially
different, and therefore, we would not obtain the continuity of the right
derivative in $\e.$ In this case, one can guarantee the existence of the
above derivatives only at the point $(x_0,\, 0)$ and continuity of the
mapping $P$ itself in a neighborhood of this point.
Then we go out of the framework of smooth problems and should use some
fine topological facts like the Brower fixpoint theorem. In this sense,
the potential of standard needle variations is rather restrictive.\,
The situation can be, seemingly, repaired if we take as $\theta$
{\it a point of approximate continuity}\, of the function $u(t),$ and
set $u_\e(t) =v$ not on the interval $(\theta -\e,\, \theta),$ but on its
intersection with that set along which $u(t)$ is continuous at the point
$\theta.$ Yet still we then obtain just the strict differentiability at
$\e=0,$ not the existence and continuity of the derivative at $\e >0.$
Moreover, this approach would require rather cumbersome constructions and
deep facts of the theory of measure, which is strongly undesirable in
a lecture course.

Nevertheless, even in the case of measurable $u(t),$ the smoothness of
the mapping $P$ can be obtained rather easily if the needle variation is
made not {\em by replacing}\, the given control $u(t)$ with a chosen
value $v$ on a time interval of length $\e$ near the point $\theta,$ but
{\em by the extension}\, of this point into an interval of length $\e$
and {\em insertion}\, the control value $v$ on this interval. (Here, the
total interval $[t_0, t_1]$ would extend to $[t_0, t_1 +\e].)$ In fact,
this would be exactly the result of the above mentioned $v-$change of time.
(The letter $v$ is involuntarily used here in two different meanings.)
\normalsize \bs

Next, we need to differentiate a scalar function of the form
$G(a,\e) = g(a,\, P(a,\e)),$ where $g(x_0,\, x_1)$ is a differentiable function
in a neighborhood of the point $(\hx_0,\, \hx_1)$ in the space $\R^{2n}.$
It is convenient to consider this as a separate property. \,
\ms

As before, let a pair  $(x(t), u(t))$ be the solution to equation (\ref{xua1}).\,
For this pair, introduce a Lipschitz continuous function $\p(t)$ (a row vector
of dimension $n)$ as the solution to the Cauchy problem
\beq{adeq}
\dot\p =\; -\p\,f_x(x,u), \qq \p(t_1) =\, -g_{x_1}(x_0,\,x_1).
\eeq
Note that $\p(t)$ does not depend on the needle variation and its
parameters $(\theta, v).$
\ms

%%------------ Лемма 2 --------------------------------------
\ble{lem2}\q
The derivative of function $G$ with respect to  $a$ and its left derivative
with respect to $\varepsilon$ at the point $(x_0,\,0)$ are given
by the formulas:
$$
\phantom{tttttttttt}  G_a(x_0,\,0)\,\bar a\, =
\,(g_{x_0}(x_0,x_1)- \p(t_0))\,\bar a,
\qq \all \bar a \in \R^n,
$$ $$
G_\e^+(x_0,\,0) =\; -\p(\theta)\cdot \D f(\theta, v).
$$
\ele

\Proof Let $\bx$ be any solution of the variational equation
$\dot{\bx} =\, f_x(x,\, u)\,\bx\,$ on some interval $[t', t''],$ no matter
with which initial value. It is easy to see that the product $\p(t)\,\bx(t)$
is constant on this interval. Indeed, its derivative
$$
\frac{d}{dt}\,(\psi\, \bx)=\; \dot\psi\,\bx\,+ \,\psi\dot{\bx}\; =\;
-\,\psi f_x\,\bx\, +\, \psi f_x\,\bx\; =0.
$$
(This fact is valid for any solutions to the linear equation
$\dot{\bx} = A(t)\,\bx$ and to its {\it adjoint}\, equation
$\dot \psi = -\psi\, A(t),$ with any integrable matrix  $A(t)$.)
\ms

Now, take an arbitrary $\bar a \in \R^n$ and let $\bx$ be the solution to
equation (\ref{urvar1}). Then, in view of Lemma  \ref{lem1},
$$
G_a(\hx_0, 0) = \; g_{x_0}\,\bar a\,+\, g_{x_1}\,P_a(\hx_0,0)\,\bar a\; =
\; g_{x_0}\,\bar a\,+\,  g_{x_1}\,\bx(t_1) \;=
$$ $$
\, = g_{x_0}\,\bar a\, -\,\p(t_1)\,\bx(t_1) \;=\;
g_{x_0}\,\bar a\, -\,\p(t_0)\,\bx(t_0) \;=\;(g_{x_0}- \p(t_0))\,\bar a.
$$
Similarly:
$$
G_\e^+(\hx_0, 0) = \; g_{x_1}\,P_\e^+(\hx_0,0)\;=\; -\p(t_1)\,\bx(t_1)\;=
$$ $$
=\; -\p(\theta)\,\bx(\theta) \;=\; -\p(\theta)\cdot \D f(\theta, v).
\qq\qq \hfill \mbox{\ctd}
$$ \ms

Thus, we examined differential properties of elementary needle variation.
However, as is well known, to prove the MP in problem (\ref{4})--(\ref{1}),
one needle variation is not enough.   Therefore, we define
\ms

{\bf The concept of package of needle variations.}\q Consider again a base
pair, which we now denote by $(\hx(t), \hu(t)),$ with a piecewise continuous
control $\hu(t)$ which, as before, is assumed to be left continuous at all
points in the semi-interval $(t_0,t_1]$ and right continuous at $t= t_0\,.$

Let $\cal N$ be an arbitrary finite set of pairs $(\theta_i,\,v_i),$ $i=1,\ldots,s,$
where $\theta_1 \le \ldots \le \theta_s$ are arbitrary points in the interval
$(t_0,t_1),$ and $v_1\,, \ldots,\, v_s\;$ are arbitrary points of the set $U.$
Let be given a vector $\varepsilon=(\varepsilon_1,\ldots,\varepsilon_s)$
with small nonnegative components $\varepsilon_i\,.$ We construct a variation
of the control that takes values $v_i$ on semi-intervals $\Delta_i$ of
length $\e_i$ near the points $\theta_i\,.$
\ssk

Define these semi-intervals as follows. If all the points $\theta_i$ are
different, we set $\D_i = (\theta_i -\e_i,\, \theta_i],$ $i=1,\ldots,s.$
If there are duplicates among $\theta_i\,,$ then for each number of duplicate
values $\theta_k = \ld = \theta_{k+p}$ we place sequentially adjacent intervals
of length $\e_{k+i}$ to the left from the point $\theta_k\,,$ i.e., we set
$$
\D_k = (\theta_k -\e_k,\, \theta_k),\qq \D_{k+1} =
(\theta_k -\e_k -\e_{k+1},\, \theta_k -\e_k), \qq \mbox {etc.},
$$
and for all non-duplicate values $\theta_i$ we set, as before,
$\D_i = (\theta_i -\e_i,\, \theta_i].$ Thus, for the given collection
$\cal N$ and the vector $\e\ge 0$ we define semi-intervals $\,\D_i$ of
length $\e_i\,.\;$ If the vector $\varepsilon$ is sufficiently small, these
semi-intervals do not overlap and all lie in the interval  $[t_0,t_1].$
Denote by $\t_i(\e)$ the right end of semi-interval $\D_i\,.$ Obviously,
it continuously depends on the vector $\e.$
\ssk

Now, define the control
$$
u_\e(t)\,=\;
\left\{\begin{array}{l}
v_i\,,\quad \mbox{ if}\quad t\in \Delta_i\,,\quad i=1,\ldots,s, \\[4pt]
\hu(t),\quad \mbox{if}\quad t\in[t_0,t_1] \setminus \bigcup_{i=1}^s\Delta_i\,.
\end{array}\right.
$$
The family of functions  $u_\e(t)$ with $\e\to 0+$ is called {\it a package
of needle variations}\, of the base control $\hu(t).$ \ssk

Substituting $u_\e(t)$ into control system (\ref{B4}), we obtain the Cauchy
problem (\ref{xea1}) with a certain initial value $a.$ As in the case of
elementary needle variation, the theory of ODE says that, for sufficiently
small $\e\ge 0$ and an initial value $a$ sufficiently close to  $\hx(t_0),$
the solution of (\ref{xea1}) exists, unique, and depends continuously on
the pair  $(a,\e).$
\ssk

Thus, we obtain the  mapping
$$
P:\; (a,\e)\; \longmapsto\; x_\e(t_1) \in \R^n,
$$
which is defined and continuous on $\W(\hx_0) \times (\calO \cap \R^s_+),$
where $\W(\hx_0)$ is a neighborhood of the point  $\hx_0 \in \R^n,$ by
$\calO$ is denoted a neighborhood of zero in  $\R^s,$ and $x_\e(t_1)$ is
the value of solution to (\ref{xea1}) at the point  $t=t_1\,.$

Lemma  \ref{lem1} implies that the mapping  $P$ is not only continuous,
but also smooth in the following sense.

%%----------- Лемма 3 --------------------
\ble{lem3} At any point $(a,\e)$ of its domain, the mapping  $P$ has the
derivative with respect to $a$ and the right derivative with respect to
each $\e_i\,,$ which continuously depend on the pair $(a,\e).$
\ele

\Proof The existence of these derivatives (and their expression through
concrete formulas) was established in Lemma  \ref{lem1}. Further, under
variations of the pair $(a,\e),$ the control $u_\e$ continuously varies
in the norm of space $L_1[t_0, t_1],$ the corresponding $x_\e$ continuously
varies in the norm of space $C[t_0, t_1],$ and for each $i-$th needle,
the pair of its "base"\, values $(\t_i(\e),\, u_\e(\t_i(\e))$ continuously
varies under variation of $\e$ by virtue of left continuity of the function
$u(t).$ Therefore, by Lemma \ref{lem1} the derivatives of mapping $P$
continuously depend on the pair $(a,\e)$ in a neighborhood of the point
$(\hx_0\,,0).$

Explain this in an example of two needles at the point $\theta$ with values
$v_1, v_2$ on semi-intervals $\D_1 = (\theta -\e_1,\, \theta], \,$
$\D_2 = (\theta -\e_1 -\e_2,\, \theta -\e_1].$ By Lemma \ref{lem1}, here
$P_{e_1}^+$ is given by the formula for the elementary needle at the point
$\theta,$ and $P_{e_2}^+(\e_1, \e_2) = \bx(t_1),$ where $\bx(t)$ is the
solution to the Cauchy problem (\ref{urbx1}) on the interval
$[\theta -\e_1 -\e_2,\,t_1]\,$ with the initial condition  \\[4pt]
\hspace*{10mm}
$ \bx(\theta -\e_1 -\e_2) =\; f(\hx(\theta -\e_1 -\e_2),\, v_2)
- f(\hx(\theta -\e_1 -\e_2),\, \hu(\theta -\e_1 -\e_2)).\,$
\\[4pt]
Both these derivatives continuously depend on $(\e_1,\e_2)$ up to the
boundary of quadrant $\R^2_+\,.$ In particular, $P_{e_2}^+(\e_1, 0)$
corresponds to the initial value
$$
\bx(\theta -\e_1) =\; f(\hx(\theta -\e_1),\, v_2) -
f(\hx(\theta -\e_1),\, \hu(\theta -\e_1)),
$$
and $P_{e_1}^+(0, \e_2)$ to the initial value
$$
\bx(\theta -\e_2) =\;
f(\hx(\theta -\e_2),\, v_2) - f(\hx(\theta -\e_2),\, \hu(\theta -\e_2)).
\qq \Box
$$
\ssk

%and their continuity follows from the fact that, under variations of vector
%$\e,$ the control $u_\e$ varies continuously in the norm of the space
%$L_1[t_0, t_1],$ and then the corresponding $x_\e$ varies continuously
%in the norm of the space  $C[t_0, t_1].$ \ctd

We will now show that if  $\e \in \inter\,\R^s_+$ (briefly, $\e>0$), then
there exist the usual partial derivatives $P_{\e_i}(a,\e)$ which coincide
with the right derivatives. This is a one-dimensional fact;\, it follows
from the below statement, which is a simple exercise in mathematical analysis.

\ble{lem4}
Let a function $\f: \R_+ \to \R$ have, for all $x\ge 0,$ the right derivative
$\f'_{\mbox{\rm \small +}}(x)$  which is continuous on  $\R_+\,.$ Then $\f$
is differentiable at all $x>0,$ and hence, $\f'(x) = \f'_{\mbox{\rm \small +}}(x).$
\ele

The proof is given in Appendix. \bs

From this lemma and Lemma \ref{lem3} it follows that the mapping $P$ is
continuously differentiable on the domain $\e>0,$ more precisely on
$\W(\hx_0) \times (\calO \cap \inter\R^s_+)),$ and its derivative has a
limit at those points $(a,\e)$ where $\e\in \pd\R^s_+\,,$ which we also
denote by $P'(a,\e).$ This, in turn, implies that $P$ is {\em strictly
differentiable}\, at the point $(\hat x_0,0),$ i.e., the difference
$P(a,\e) - P'(\hat x_0,0)(a,\e)$ is Lipschitz continuous in a neighborhood
of $(\hat x_0,0)$ (of course, intersected with $\R^n \times \R^s_+),$
with the constant tending to zero together with the radius of the neighborhood.
(The concept of strict differentiability, we believe to be presented earlier
in this course of optimal control.) \ssk

Thus, the smoothness of mapping $P,$ i.e. the continuity of its partial
derivatives, is reduced to the continuity of the derivatives of one elementary
needle in dependence of the base process. We borrowed this idea from
Magaril--Il'yaev \cite{Mag}.
\ms

The key point of our proof is that the mapping $P$ can be extended to a
"full"\, neighborhood of  $(\hat x_0,0),$ i.e., to arbitrary $a \in \W(\hx_0)$
and arbitrarily small values of the vector $\e,$ preserving its strict
differentiability at $(\hat x_0,0).$ This is due to the following general
statement. \ssk

%%----------------------------------------------------------------------
\ble{lem5}
Let $K$ be a closed convex cone in  $\R^s$ with a nonempty interior, and
a mapping $P: \R^n \times K \to \R^m$ be strictly  differentiable at $(0,0).\,$
Then this mapping can be extended to a mapping $\wt P: \R^n \times \R^s \to \R^m$
still strictly differentiable at  $(0,0).\,$
\ele

The proof is carried to Appendix. (As the cone $K$ we have the nonnegative
orthant $\R^s_+$.)\, The obtained extension will be denoted by the same
letter $P.\,$ Note that for  $K = \R^s_+$ and smooth $P$ (like in our case),
one can do without the concept of strict differentiability, using the
standard concept of continuous differentiability (see \cite{Kor}):
\ssk

\ble{lem5a}
Let a mapping $P: \R^n \times \R^s_+ \to \R^m$ be continuously differentiable
in the interior of its domain, and its derivative be continuous up to the boundary.\,
Then $P$ can be extended to a mapping $\wt P: \R^n \times \R^s \to \R^m\,$
continuously differentiable in the whole space.
\ele

The proof is also carried to Appendix. \ssk

We draw attention once again that all the previous arguments concern only
differential properties of the control system on the package of needle
variations and are not directly related to optimization problem. So they
are of inherent interest and can be used in other areas of control theory.
\ms

Now we are ready to prove the MP\, for problem B.
\bs

%%-------------------------------
\subsection*{Maximum principle for the package of needle variations}

Let a process $\hw = (\hx, \hu)$ provide the strong minimum in problem B.\,
Fix some package  ${\cal N}$ of needle variations. Using lemma \ref{lem5}
(or lemma \ref{lem5a}) extend the mapping $P$ to a full neighborhood of
the point $(\hx_0,0),$ where $\hx_0 = \hx(t_0),$ preserving its strict
differentiability at this point (or, respectively, its continuous
differentiability in a neighborhood of this point). \ssk

In the finite dimensional space  $\R^n\times \R^s$ with elements  $(a,\e)$
consider the following  {\em problem\,} $Z(\calN)$:
$$
J(a,P(a,\varepsilon)) \to \min,  \quad F(a,P(a,\varepsilon))\le 0,
\quad K(a,P(a,\varepsilon))=0, \quad -\varepsilon\le 0.
$$
This problem is a restriction of problem B, and therefore the fact that the pair
$(\hat x,\hat u)$ provides the strong minimum in problem B\, implies that the pair
$(\hat a= \hx_0,\,\,\hat\e=0)$ provides the local minimum in problem $Z(\calN)$
and hence satisfies the stationarity conditions.\,
\ssk

Let us write down these conditions for the point $(\hx_0,0)$ in
problem $Z(\calN).$ They say that there exist  Lagrange multipliers
$$
\alpha_0\in\R,\quad \alpha\in {\R^{d(F)}},\quad \beta\in {\R^{d(K)}},
\quad \g\in\R^s,
$$
such that the following conditions hold
\begin{eqnarray} && \label{stat1a}
\alpha_0\ge 0,\quad \alpha \ge 0, \quad \g \ge 0,\\[4pt]
&& \label{stat2a} \alpha_0+|\alpha|+ |\beta|+ |\g| >0,\quad
\alpha F(\hat x_0,\hat x_1)=0,\\[4pt]
&& \label{stat3a}{\cal L}_{a}(\hat x_0,0)=0,\qquad
{\cal L}_\varepsilon(\hat x_0,0)=0,
\end{eqnarray}
where
$$
{\cal L}(a,\e)= \big(\alpha_0 J + \alpha F +\beta K\big)(a,P(a,\e))\, -\,\g\,\e\;
= \;  { l} (a,P(a,\e))\,- \g\,\e
$$
is the Lagrange function for problem  $Z(\calN).$ \ssk

Note once again, that unlike in some other proofs of the MP, including
\cite{Pont, Bolt, Mag}, we do not construct the cone of variations $\bx(t_1),$
but simply consider the finite-dimensional problem correspon\-ding to
the package of needle variations, and use already known (and presented
earlier in the course) stationarity conditions for this problem. \ms

Now, decipher  conditions (\ref{stat3a}).\, Let a Lipschitz continuous function
$\p(t)$ be the solution of the {\it adjoint equation}\, (to equation (\ref{urvar1})
along the optimal process $(\hx,\hu)$)
\begin{equation}\label{adj}
\dot \psi =\, -\psi f_x(\hat x(t),\hat u(t))
\end{equation}
with the boundary condition $\psi(t_1)= -l_{x_1}(\hat x_0,\hat x_1)\,.$
Note again that $\p$ depends only on the tuple of multipliers $(\a_0,\a,\b)$
(since the endpoint function $l$ is expressed through them) and does not
depend explicitly on the chosen package of variations.\, Moreover, $\p$
is uniquely determined by $(\a_0,\a,\b).$
\ssk

By Lemma \ref{lem2} applied to the function $g= \calL,$  the first equality in
(\ref{stat3a}) means that   $\all \bar a \in \R^n$
$$
{\cal L}_{a}(\hat x_0,0)\,\bar a\, = \;
(l_{x_0}(\hx_0, \hx_1)\, -\,\p(t_0))\,{\bar a}\;=\,0,
$$
which implies that  $\p(t_0) = l_{x_0}(\hx_0, \hx_1).$ Thus, the function
$\p(t)$ satisfies the boundary conditions at both ends of the interval:
\beq{tr1}
\p(t_0) = l_{x_0}\,, \qq \p(t_1) = -l_{x_1}\,.
\eeq
They are called {\it transversality conditions}.
\ms

The second equation in (\ref{stat3a}) means that for every  $i$
\beq{lei}
{\cal L}_{\varepsilon_i}(\hx_0,0)= \;
-\p(\theta_i)\cdot \D f(\theta_i, v_i) - \g_i\; = 0,
\eeq
i.e.,  $\q \p(\theta_i)\cdot \D f(\theta_i, v_i)\, =\;
\p(\theta_i) \Bigl( f(\hx(\theta_i),v_i)) -
 f(\hx(\theta_i), \hu(\theta_i)) \Bigr)\, = - \g_i \lee 0.$
\bs

Introduce the Pontryagin function  $H(\p,x,u) = \p\,f(x,u).$ Then the
last condition means that the following {\em finite valued maximum condition}\,
holds:
\beq{usmax}
H(\p(\theta_i),\hx(\theta_i),v_i) \lee
H(\p(\theta_i),\hx(\theta_i), \hu(\theta_i)),
\qq  i=1, \ld, s,
\eeq
and the adjoint equation (\ref{adj}) can be written as
\beq{adjh}
\dot\p(t) =\, - H_x (\p(t),\hx(t), \hu(t)).
\eeq

(Condition (\ref{usmax}) can be interpreted as follows:\, for any $\theta_i$
from the given package, the function $H(\p(\theta_i),\hx(\theta_i),u)$ takes
its maximum over all values of the control $u$ presenting in this package
at the point $\theta_i\,,$ at the optimal $u = \hu(\theta_i).$)
\ssk

Thus, for any given package $\calN,$ we obtain a tuple of  Lagrange multipliers
which generates a function  $\p(t)$ such that  conditions (\ref{stat1a}),
(\ref{stat2a}), (\ref{adj}), (\ref{tr1}), (\ref{usmax}) hold.\,  This tuple
of Lagrange multipliers, in general, depends on the package.  Conditions
(\ref{stat1a}), (\ref{stat2a}), (\ref{adj}), (\ref{tr1}) are the same for
all packages, while condition (\ref{usmax}) is directly related to the given
package. Our goal now is to pass to the maximum condition for all
$t \in (t_0, t_1)$ and $v \in U$ with a tuple of multipliers independent
of $t$ and $v$.
\ms

Note preliminarily that, if $\alpha_0+|\alpha|+|\beta|=0,$ then $l=0,$
hence $\p(t_1)=0,$ and since $\p$ satisfies a homogeneous linear equation,
$\p(t) \equiv 0,$ and then by (\ref{lei}) we obtain $\g =0,$ which contradicts
condition (\ref{stat2a}).\, Therefore, everything is determined by the
multipliers $(\a_0, \a, \b),$ and the nontriviality condition can be written
as $\alpha_0 +|\alpha| +|\beta| >0.$ We replace it by the normalization
condition $\,\alpha_0 +|\alpha| +|\beta| =1.$
\bs

%%-----------------------
{\bf Arrangement of optimality conditions, passage to the universal MP.}\\[4pt]
Denote by  $\Lambda({\cal N})$ the set of all tuples of multipliers
$(\alpha_0,\alpha,\beta)$ satisfying the conditions
\begin{eqnarray} && \label{stat1aa}
\alpha_0\ge 0, \quad \alpha\ge 0,\quad \alpha_0+|\alpha|+|\beta|=1,
\quad \alpha F(\hat x_0,\hat x_1)=0, \\[4pt]
&& \label{stat2aa}
\dot \psi =\, -\psi f_x(\hat x(t),\hat u(t)), \qq
\p(t_0) = l_{x_0}\,, \q \p(t_1) = -l_{x_1}\,,
\end{eqnarray}
and the finite-valued maximality condition (\ref{usmax}).\, (Such a tuple
may be not unique.)\,\\ It is easy to see that $\Lambda({\cal N})$ is a
nonempty compact set in the space  $\R^{1+d(F)+ d(K)}.$ \ssk

Thus, taking all possible packets $\calN,$ we obtain for each of them a
nonempty compact set $\Lambda({\cal N}).\,$ Let us show that the family
of all these compact sets $\{ \Lambda({\cal N})\}_{\cal N}$ is {\em centered},\,
i.e. has {\em the finite intersection property}.\,
To this end, we introduce an ordering in the set of all packages. We say
that ${\cal N}_1\subset{\cal N}_2\,$ if each pair $(\theta_i,v_i)$ from
${\cal N}_1$ belongs (possibly with a different index $i$) also to
${\cal N}_2\,.$ It is clear that for any two packages ${\cal N}_1$ and
${\cal N}_2\,,$ there exists a third one containing each of them, e.g.,
their union.  Further, it is clear that the expansion of $\calN$ narrows
the set $\L(\calN),$ i.e., the inclusion ${\cal N}_1\subset{\cal N}_2$
entails the reverse inclusion $\L({\cal N}_1)\supset \L({\cal N}_2).$
Now, let be given a finite number of compacta
$\Lambda({\cal N}_1),\ldots,\Lambda({\cal N}_r).$ Take any package
${\cal N}$ containing all the packages ${\cal N}_1,\ldots,{\cal N}_r\,.$
Then the nonempty compact set $\Lambda({\cal N})$ is contained in each
of the sets $\Lambda({\cal N}_1),\ldots,\Lambda({\cal N}_r)$ and hence
in their intersection. This implies the finite intersection property
of the system $\{ \Lambda({\cal N})\}_{\cal N}\,$ and hence the nonemptyness
of its total intersection
$$
\L_*\; =\; \bigcap\limits_{{\cal N}}\,\Lambda({\cal N}).
$$

Take an arbitrary tuple of multipliers  $(\alpha_0,\alpha,\beta)\in \L_*$
and let $\psi(t)$ be the adjoint variable corresponding to this tuple.
%where the gradient $l_{x_1}(\hat x_0,\hat x_1)$ corresponds to the tuple
%$(\alpha_0,\alpha,\beta).$
By definition, this tuple satisfies conditions (\ref{stat1aa})--(\ref{stat2aa}),
while maximality condition (\ref{usmax}) for the Pontryagin function
with respect to the control holds now {\it for any pair}
$(\theta,v) \in (t_0, t_1) \times U,$
since any such pair belongs to some package.  Since the function $u(t)$
is piecewise continuous, this condition as well holds for the endpoints
of the interval.

Thus, conditions $(i') - (v')$ of theorem \ref{th2}\, are proved.
\bs

%%-----------------------------------
{\bf Constancy of the Pontryagin function along the optimal process.}\\
It remains to prove that along the optimal process
\begin{equation}\label{const-h}
H(\psi(t),\hat x(t),\hat u(t))=\, \const \quad\mbox{ a.e.\, on }\; [t_0,t_1],
\end{equation}
i.e. condition $(vi')$ is satisfied. Then condition $(vii')$ would follow
from the maximality condition (\ref{usmax}).\, Let us show that condition
(\ref{const-h}) follows from already proved condi\-tions of MP.\,
This can be done in two ways:\, 1) to show this directly, and\,
2) to get (\ref{const-h}) automatically from the other conditions of MP\,
after a passage to a problem on a nonfixed time interval.

In our opinion, the second way is more simple and natural (see the next
section).\, A merit of the first way is that it establishes a direct
dependence of condition (\ref{const-h}) from other conditions of MP.\,
This way was proposed in the books \cite{Pont, Bolt}.\, Here we give a
modification of considerations in \cite{Bolt} (those in \cite{Pont} are
more complicated). \ms

Consider the function $h(t,v) = H(\psi(t),\hat x(t),v) = \p(t)\,f(\hat x(t),v).$
We have to show that $h(t,\hu(t))$ is constant on $[t_0,t_1].$
Let us argue as follows.
\ssk

First of all, by the maximality condition,
\beq{h1}
h(t,v) \lee h(t,\hu(t)) \qq \all t\in [t_0,t_1], \q \all v\in U.
\eeq

Further, the function $h(t,v)$ and its partial derivative $h_t(t,v)$
are continuous on $[t_0,t_1]\times U,$ and at any point of continuity
of the control $\hat u(t),$ we have
\beq{h2}
h_t(t,\hu(t))\;=\; 0.
\eeq
Indeed,
$$
h_t(t,\hu(t)) =\; \left.\Big(\dot\p(t)\,f(\hat x(t),v)\, +\,
\p(t)\,f_x(\hat x(t),v)\,\dot{\hx}(t)\Big)\right|_{v= \hu(t)} \; =
$$ $$
=\; \left.\Big(-\p(t)\,f_x(\hat x(t),\hu(t))\,f(\hat x(t),v)\,+\,
\p(t)\,f_x(\hat x(t),v)\,f(\hat x(t),\hu(t))\Big)\right|_{v=\hu(t)}\; =0.
$$

Let us now show that, at any point of continuity of $\hat u(t),$
the total derivative $\frac d{dt}\, h(t,\hu(t))$ also equals zero.
Take any such point $t,$ and let two points $\t_2 > \t_1$ converge to
$t.$ Denote $\D\t=\tau_2-\tau_1\,.$ We have
$$
h(\t_2, \hu(\t_2)) - h(\t_1,\hu(\t_1)) =\;\;
h(\t_2, \hu(\t_2))- h(\t_1, \hu(\t_2))  \;+ \phantom{pppppppp}
$$ \beq{ht12}
\phantom{ppppprrrrpyyy} +\; h(\t_1, \hu(\t_2)) - h(\t_1, \hu(\t_1)).
\eeq
The modulus of the first difference in the right hand side of this equality
can be estimated by the mean value theorem:
$$
\Big|h(\t_2, \hu(\t_2))- h(\t_1, \hu(\t_2))\Big| \lee
\max_{\t\in [\t_1, \t_2]} \Big| h_t(\t,\hu(\t_2)) \Big|\cdot|\D\t| \; =\; o(\D\t),
$$
since the above maximum tends to $|\, h_t(t,\hu(t))| =0.\,$ The second
difference in (\ref{ht12}) is nonpositive by virtue of maximality condition
(\ref{h1}). Thus, $h(\t_2, \hu(\t_2)) - h(\t_1,\hu(\t_1))\le o(\D\t).$
\ssk

On the other hand, the same difference admits another representation:
$$
h(\t_2, \hu(\t_2)) - h(\t_1,\hu(\t_1)) =\;\;
h(\t_2, \hu(\t_2))- h(\t_2, \hu(\t_1))  \;+ \phantom{pppppppp}
$$ \beq{ht21}
\phantom{ppppprrrrpyyy} +\; h(\t_2, \hu(\t_1)) - h(\t_1, \hu(\t_1)).
\eeq
from which, by the same reasons, we obtain
$\; h(\t_2, \hu(\t_2)) - h(\t_1,\hu(\t_1))\ge o(\D\t).$ \\
Hence,
$$
h(\t_2, \hu(\t_2)) - h(\t_1,\hu(\t_1))\; =\; o(\D\t).
$$
This implies that the total derivative $\dot h(t,\hu(t))$ exists and
equals zero at all points of continuity of control $\hat u(t),$ i.e.
$h(t,\hu(t)) = \const$ on any interval of continuity of $\hat u(t).$
Let us show that this constant is one and the same for all intervals.
\ssk

Let $\tau\in (t_0,t_1)$ be a discontinuity point of control $\hat u(t)$
(according to our assumption, there is no more than a finite number of
such points).\, Show that $h(t,\hu(t))$ does not have a jump at this point.\,
For $t>\tau,$ by the condition (\ref{h1}), we have
$$
h(t,\hu(t)) \gee h(t,\hu(\t-0)).
$$
Passing to the limit as  $t \to \tau+0\,,$ we obtain
$$
h(\t,\hu(\t+0)) \gee h(\t,\hu(\t-0)).
$$
Similarly, for $t<\tau,$ by the condition (\ref{h1}), we have
$$
h(t,\hu(t)) \gee h(t,\hu(\t+0)).
$$
Passing to the limit as $t\to\tau-0\,,$ we obtain
$$
h(\t,\hu(\t-0)) \gee h(\t,\hu(\t+0)).
$$
Therefore, $h(\t,\hu(\t-0)) = h(\t,\hu(\t+0)).$
This and the preceding imply that $h(t,\hu(t))= \const\,$ on the entire
interval $[t_0, t_1].$  Thus, condition (\ref{const-h}) is proved, and
with it, theorem \ref{th2}\, is completely proved.
\bs

%%------------------------------------------------
\section{Proof of the maximum principle for the problem\\
on a nonfixed time interval}

Now we get back to problem (\ref{4})--(\ref{1}) on a nonfixed time interval,
called problem A.\, Let $w^0=(x^0(t), u^0(t))\mid  t\in [\hat t_0, \hat t_1])$
be an admissible process in this problem which provides the strong minimum.
(Here, it will be convenient to equip $x,\,u$ with the superscript $0.$)\,
As before, we assume the control $ u^0(t)$ to be piecewise continuous.\,
\ms

Let us reduce problem A\, to an accessory problem on a {\em fixed}\, time
interval, in order to use then Theorem~\ref{th2}.  To this end, we introduce
a "new time" $\tau$ varying in the fixed interval
$[\tau_0, \tau_1]:=[\hat t_0,\hat t_1],$ while the original time $t$ will
be regarded as one more state variable $t= t(\t)$ satisfying the equation
\beq{dttau}
\frac {dt}{d\t}\; =\; v(\t),
\eeq
where $v(\t)>0$ is one more control variable. The original variables $x,\, u$
in the new time are $\tl x(\t) = x(t(\t)),$ $\tl u(\t)= u(t(\t));\,$
they satisfy the equation
\beq{eq1tau}
\frac {d\tl x}{d\t}\; =\; v(\t)\,f(t(\t), \tl x(\t), \tl u(\t)).
\eeq
The endpoint block of the problem now is:
\beq{eq2tau}
K(t(\t_0),\tl x(\t_0), t(\t_1), \tl x(\t_1)) =0, %\qq j=1,\ld, k,
\eeq
\beq{eq3tau}
F(t(\t_0),\tl x(\t_0), t(\t_1), \tl x(\t_1)) \le 0, %\qq i=1, \ld, \nu,
\eeq
\beq{eq4tau}
J = F_0(t(\t_0),\tl x(\t_0), t(\t_1), \tl x(\t_1)) \to \min.
\eeq \\[1pt]
We see that it depends only on the endpoint values of the state variables
$t$ and $\tl x.$ Problem (\ref{dttau})--(\ref{eq4tau}) will be called
problem $\wt A.$ Here $t(\tau),\; \wt x(\tau)$ are the state variables,
and $v(\tau),\; \wt u(\tau)$ are the controls. The condition $v>0$ should
be regarded as a constraint on the new control. \ssk

With the process  $w^0$  in problem  A\, we associate the process
$$
\wt w^0\,=\;( v^0(\tau),\, t^0(\tau),\, \wt x^0(\tau),\, \wt u^0(\tau))\;
\mid \; \tau\in [\tau_0, \tau_1])
$$
in problem  $\wt A,\;$ where  \vad
\begin{equation}\label{3nn}
v^0(\tau)\equiv 1, \q t^0(\tau)\equiv\tau, \quad
\wt x^0(\tau)\equiv x^0(\tau), \quad
\wt u^0(\tau)\equiv u^0(\tau).
\end{equation}
Obviously, this process is admissible in problem $\wt A.$ Moreover, the
fact that process $w^0$ provides the strong minimum in Problem $A$ entails
that process $\wt{w}^0$ provides the strong minimum in Problem  $\wt A.$
It is easily obtained by contradiction. By Theorem~\ref{th2}, $\wt {w}^0$
satisfies the MP.\, Let us write down its conditions for process $\wt {w}^0$
in Problem $\wt A,$ omitting now the tilde and taking into account that
$dt^0/ d\t \equiv 1$ and $t^0(\t) \equiv \t,$ so that the derivative with
respect to $\t$ can still be denoted by the upper dot.
\ssk

The endpoint Lagrange function in Problem $\wt A$ is the same as in
Problem~$A$:  \vad
$$
l=\; \alpha_0 F_0+\alpha F+\beta K,
$$
while the Pontryagin function in Problem  $\wt A$ has the form
$$
\wt H(\p_x,\p_t,\,t,x,u,v) =\; (\p_x\,f(t,x,u) + \p_t)\,v.
$$

The conditions of maximum principle for the process $\wt{w}^0$ in Problem
$\wt A$ are as follows: \, there exists a number  $\alpha_0\,,$ vectors
$\alpha \in\R^{d(F)},$ $\beta \in \R^{d(K)},$ Lipschitz continuous functions
$\psi_t$ and $\psi_x$ of dimension 1 and $n,$ respectively, such that
the following conditions hold: \ssk

\noi nonnegativity and nontriviality  \vad
\begin{equation}\label{mp1}
\alpha_0\ge0,\quad \alpha\ge0;\quad  \alpha_0+|\alpha|+|\beta|>0,
\eeq
complimentary slackness   \vad
\beq{dopnez}
\alpha F(t^0(\tau_0), x^0(\tau_0), t^0(\tau_1), x^0(\tau_1))=0;
\end{equation}
adjoint equations   \vad
\beq{dpsixtau}
-\dot\psi_x(\t) = \wt H_x \;=\; \p_x(\t)\,f_x(t^0(\t),x^0(\t),u^0(\t)),
\eeq
\beq{dpsittau}
-\dot\psi_t(\t) = \wt H_t \;=\; \p_x(\t)\,f_t (t^0(\t),x^0(\t),u^0(\t)),
\eeq
transversality
\beq{transxtau}
\p_x(\t_0) =\; l_{x(\t_0)}\,, \qq  \p_x(\t_1) =\; -l_{x(\t_1)}\;,
\eeq
\beq{transttau}
\p_t(\t_0) =\; l_{t(\t_0)}\,, \qq  \p_t(\t_1) =\; -l_{t(\t_1)}\;,
\eeq  \\
and maximality with respect to $u$ and $v$ \vad
\begin{equation}\label{mp10}
v\Big(\psi_x(\tau)\,f(t^0(\tau), x^0(\tau),u)+ \psi_t(\tau)\Big)\le
v^0(\tau)\Big(\psi_x(\tau)\,f(t^0(\tau), x^0(\tau), u^0(\tau))+\psi_t(\tau)\Big)
\end{equation}
$$
\mbox{for all}\q u\in U, \q v>0,\;\q \tau\in[\tau_0,\tau_1].
$$ \ssk

\noi
For $u = u^0(\t),$ the maximality of $\wt H$ over $v>0$ at the point $v^0 =1$
implies that
\begin{equation}\label{mp13}
\psi_x(\tau)\,f(t^0(\tau), x^0(\tau), u^0(\tau))+\psi_t(\tau)\,=\, 0\quad
\mbox{a.e.\, on}\quad [\tau_0,\tau_1],
\end{equation}
and then (\ref{mp10}) means that
\begin{equation}\label{mp14}
\psi_x(\tau)\,f(t^0(\tau), x^0(\tau), u^0(\tau)) +\psi_t(\tau) \le 0\quad
\mbox{for all}\;\; u\in U, \;\; \tau\in[\tau_0,\tau_1].
\end{equation}
Taking into account that  $t^0(\tau) =\t,$ we see that all conditions
$(i)-(vii)$ of the maximum principle for the process $w^0$ in Problem $A$
are obtained.\,  Theorem~\ref{th1}\, is completely proved.  \ctd \bs

Note that, for the autonomous problem, relation (\ref{dpsittau}) implies
$\psi_t(\tau)=\const,$ whence (\ref{mp13}) means that
$H(\p_x(t), x^0(t), u^0(t)) = -\p_t = \const,$ hence, as was said above,
we automatically obtain condition (\ref{const-h}).
\bs

%%-------------------------------

\section{Appendix}

To prove Lemma  \ref{lem2},\, we need the following
\ble{lem7}
Let $f$ be a continuous function on an interval $[a,b],$  such that
for all $x \in [a,b)$ there exists the right derivative
$f'_{\mbox{\rm \small +}}(x) >0.\,$ Then $f(b) > f(a).$
This also implies that $f$ strictly increases on $[a,b].$
\ele

\Proof Suppose that $f(a)=0.$ Since $f'_{+}(a) >0,$ we have $f(x)>0$ on
the interval $(a,\,a+\d)$ for some $\d>0.$ Take any point $a_1$ of this
interval, and let $a_2$ be the maximal of the points $x$ for which
$f(x) \ge f(a_1)$ on the interval  $[a_1,\,x].$ If $a_2 <b,$ then, since
$f'_{\mbox{\rm \small +}}(a_2) >0,$ we obtain $f(x) > f(a_2) \ge f(a_1)$
in some right half-neighborhood of $a_2\,.$ This contradicts the maximality
of $a_2\,.$ Therefore, $a_2 =b,$ and hence the lemma is proved. \ctd
\bs

Using this lemma, we can prove an analogue of Lagrange mean value theorem
for the right derivative.

\ble{lem8}
Let $f$ be a continuous function on an interval $[a,b].$ Assume that
for all $x \in [a,b)$ there exists the right derivative
$f'_{\mbox{\rm \small +}}(x)$ which is continuous in $x.\,$ Then
$$
f(b) - f(a)\, =\, f'_{\mbox{\rm \small +}}(\tl x)\,(b-a) \qq
\mbox{for some}\q \tl x\in [a,b).
$$
\ele

\Proof  Subtracting a linear function from  $f,$ we come to the case
$f(a) = f(b) =0$ (i.e., we have to prove an analogue of Rolle's theorem).
Assume the assertion of lemma is false;\, say,
$f'_{\mbox{\rm \small +}}(x)>0$ for all $x \in [a,b).$ Then, by Lemma
\ref{lem7}, we have $f(b) > f(a),$ which contradicts the hypothesis.
\ctd \bs

{\bf Proof of Lemma \ref{lem4}}.\, Take any point  $x >0.$
Without loss of generality assume that $\f'_{\mbox{\rm \small +}}(x) =0.$
We need to show that the left derivative $\f'_{{\small -}}(x) =0,$
i.e., that $|\,\f(x) - \f(x-\d)| = o(\d)$ as $\d \to 0+.$
By Lemma \ref{lem8}}, we have
$$
|\,\f(x) - \f(x-\d)|\; =\; \f'_{\mbox{\rm \small +}}(\wt x)\,\d \qq
\mbox{for some } \q \wt x\in [x,\,x-\d).
$$
By the continuity of right derivative,
$\f'_{\mbox{\rm \small +}}(\wt x) \to \f'_{\mbox{\rm \small +}}(x) =0,$
which implies the required. \ctd \ms

{\bf Proof of Lemma \ref{lem5}}.\, Without loss of generality assume that
$P'(0,0)=0,$ i.e., the mapping $P$ is Lipschitz continuous in any
$r-$neighborhood of the point $(0,0) \in \R^n \times K$ with a constant
$\mu(r)\to 0$ as $r \to 0.$ Fix any vector  $h \in {\inter}\,K,$ and for
every $y \in \R^s$ set $\l(y) = \min\,\{\l \ge 0\;|\; y +\l h \in K\}.$
Define a mapping $\f: \R^s \to K$ by the formula  $\f(y) = y + \l(y)\,h.$
(It is the projection of $\R^s$ to $K$ along the vector $h.)$ For all
$y \in K$ we obviously have $\f(y) =y.$ Clearly, $\f$ is Lipschitz continuous
with some constant $L.$ Then the mapping $\wt P(x,y) = P(x,\,\f(y))$ is
defined on the whole  $\R^n \times \R^s$ and Lipschitz continuous in the
$r$-neighborhood of the point $(0,0)$ with the constant $L\,\mu(r)\to 0$
as $r\to 0.$ Therefore, it has a strict derivative at zero equal to zero.
\ctd
\ms

Note that this lemma is as well true in the case of a more general mapping
$P: X \times K \to Z,$ where $X,\, Y,\, Z\;$ are arbitrary normed spaces
and $K \subset Y$ is a closed convex cone with a nonempty interior.\,
The proof remains unchanged.
\ms

To prove Lemma \ref{lem5a}, let us establish, like in \cite{Kor}, the
following, in fact one-dimensional (and certainly well-known, see e.g.
\cite{Lich}), property.

\ble{lem9}
Suppose that a mapping $f: \R^n \times \R_+ \times \R^k_+ \to \R^m$ has
the derivative at all interior points of the domain, which is continuous
up to the boundary of this domain. Extend $f$ to the whole space
$\R^n \times \R \times \R^k_+\,$ by setting
$$
\wt f(x,y,z)\, =\; -f(x, -y,z)\,+\, 2 f(x,0,z) \qq \mbox{for }\;\; y<0.
$$
Then, the resulting mapping  $\wt f$ has the same smoothness, now
in the whole its domain.
\ele

The proof is an elementary check.
\ms

{\bf Proof of Lemma \ref{lem5a}}\, is a consecutive application of
Lemma~\ref{lem9}\, to all components of the vector  $\e.$
\bs

{\bf Acknowledgments.}\, The authors thank\, V.A. Dykhta\, for
valuable remarks\, and G.G. Magaril-Il'yaev for useful discussions.
\ms

%%-------------------------------------------------------------------

%\end{document}

%%=======================================

\vspace{10mm}

{\bf Andrei V. Dmitruk} -- Russian Academy of Sciences, Central Economics and
Mathematics Institute, Moscow, Russia;\, and Lomonosov Moscow State University,
Moscow, Russia.\q {\tt dmitruk@member.ams.org.} \\

{\bf Nikolai P. Osmolovskii} -- University of Technology and Humanities in Radom,
Poland; Systems Research Institute, Polish Academy of Sciences Warszawa,
Poland;  and Moscow State University of Civil Engineering, Moscow, Russia.\\
\, {\tt osmolovski@uph.edu.pl.}

\end{document}

%\newenvironment{proof}{Proof}[section]
%\def\theequation{\arabic{section}.\arabic{equation}}

%-------------------------------------------------------------------
\def\a{\alpha}   \def\g{\gamma}  \def\b{\beta}  %\def\th{\theta}
\def\l{\lambda}  \def\w{\omega}   \def\W{\Omega} \def\e{\varepsilon}
\def\d{\delta}   \def\f{\varphi}  \def\D{\Delta}    \def\r{\rho}
\def\s{\sigma}   \def\G{\Gamma}   \def\L{\Lambda}
\def\pd{\partial} \def\ex{\exists\,}  \def\all{\forall\,}
\def\dy{\dot y}   \def\dx{\dot x}   \def\t{\tau}
\def\bx{\bar x}   \def\bu{\bar u}   \def\bw{\bar w}  \def\bz{\bar z}
\def\p{\psi} \def\bh{\bar h}
\def\calO{{\cal O}}   \def\calP{{\cal P}}  \def\calR{{\cal R}}
\def\calN{{\cal N}}   \def\calL{{\cal L}}   \def\calQ{{\cal Q}}
\def\hx{\hat x}  \def\hw{\hat w}   \def\hu{\hat u}  \def\hv{\hat v}
\def\hz{\hat z}   % \def\ht{\hat t}

\def\empty{\O}     %\def\empty{\emptyset}
\def\lee{\;\le\;}   \def\gee{\;\ge\;}
\def\le{\leqslant}  \def\ge {\geqslant}
\def\Proof{{\bf Доказательство.}\quad}

\def\ov{\overline}  \def\ld{\ldots}
\def\tl{\tilde}   \def\wt{\widetilde}  \def\wh{\widehat}
\def\lan{\langle}  \def\ran{\rangle}

\def\hsp{\hspace*{8mm}} \def\vs{\vskip 8mm}
\def\bs{\bigskip}     \def\ms{\medskip}   \def\ssk{\smallskip}
\def\q{\quad}  \def\qq{\qquad}
\def\np{\newpage}   \def\noi{\noindent}
\def\pol{\frac12\,}    \def\const{\mbox{const}\,}
\def\vad{\vadjust{\kern-5pt}}
\def\R{\mathbb{R}}   \def\ctd{\hfill $\Box$}  \def\inter{{\rm int\,}}

\def\dis{\displaystyle}
\def\lra{\longrightarrow}  \def\iff{\Longleftrightarrow}
\def\beq{\begin{equation}\label}  \def\eeq{\end{equation}}
\def\bth{\begin{theorem}\label}   \def\eth{\end{theorem}}
\def\ble{\begin{lemma}\label}     \def\ele{\end{lemma}}
\def\begar{\begin{array}} \def\endar{\end{array}}

\newcommand{\vmin}{\mathop{\rm vraimin}}
\newcommand{\vmax}{\mathop{\rm vraimax}}

\newcommand{\blue}[1]{\textcolor{blue}{#1}}
\newcommand{\red}[1]{\textcolor{red}{#1}}
\newcommand{\gr}[1]{\textcolor{green}{#1}}

\def\bls{\baselineskip}   \def\nbls{\normalbaselineskip}
%\bls =1.1\nbls
%----------------------------------------------------------------

\hfill УДК 517.977.52

\begin{center}
{\large\bf О  доказательстве принципа максимума Понтрягина \\
с помощью игольчатых вариаций}
\\[8mm]
{\large А.В. Дмитрук, \quad Н.П. Осмоловский}
%{\large А.В. Дмитрук\footnotemark[1], \quad Н.П. Осмоловский \footnotemark[2]}
%\\[20pt]
%{\large Moscow State University}\\[4pt]
%\footnotemark[1]{avdmi@cemi.rssi.ru}\,\,
%\footnotemark[2]{osmolovsky@uph.edu.pl}\\[20pt]
\end{center}
\vs

{\bf Аннотация}. {\small
Предлагается доказательство принципа максимума для общей задачи оптимального
управления понтрягинского типа с помощью игольчатых вариаций путем сведения
ее к семейству гладких конечномерных задач, аргументом которых являются ширины
иголок в каждом пакете, применения для каждой задачи стандартного правила
множителей Лагранжа, и затем "спрессовывания"\, полученных необходимых
условий в одно универсальное условие оптимальности с помощью понятия
центрированной системы компактов.}
\ms

{\bf Ключевые слова:}\, {\small принцип максимума Понтрягина, пакет
игольчатых вариаций, оператор концевого значения, гладкая конечномерная
задача, множители Лагранжа, конечнозначное условие максимума,
центрированная система компактов.}
\bs

\section{Введение}
Как известно, первоначальное доказательство принципа максимума (ПМ) как необходимого
условия сильного минимума, данное в книге Понтрягина и др. \cite{Pont},
было весьма сложным и не укладывалось в рамки классических методов вариационного
исчисления. Позднее усилия многих математиков были направлены на то, чтобы
переосмыслить это доказательство и распространить ПМ на новые классы задач.
Это было сделано в целом ряде отечественных и зарубежных работ;\, см., например,
\cite{DM65} -- \cite{MDO} (разумеется, этот список ни в коей мере не претендует
на полноту).\, С другой стороны, вопрос о наиболее простом и прозрачном
доказательстве ПМ для {\it классической задачи понтрягинского типа}\, также
был предметом внимания специалистов. Это особенно важно с учебно-методической
стороны в целях преподавания курса оптимального управления как в математических,
так и в инженерных вузах. \ssk

Основная идея доказательства ПМ состоит во введении некоторого класса
(семейства) вариаций управления, по которому {\it можно дифференцировать.}
Для управляемых систем ОДУ известны следующие классы вариаций:\,
а) равномерно малые вариации,\, б) игольчатые вариации и их "пакеты"\,,\,
в) $v- $замена времени,\,  г) скользящие режимы.  \ssk

Равномерно малые вариации соответствуют {\it слабому}\, минимуму и приводят к
{\it условию стационарности,} т.е. к  уравнению Эйлера--Лагранжа.\,
Это очень важный результат, который может служить отправной позицией для
дальнейших исследований данного процесса на оптимальность, однако он далек
от ПМ. Для получения условий {\it сильного}\, минимума этих вариаций явно
недостаточно.  \ssk

Простейшим вариантом "немалой"\, вариации является {\it игольчатая вариация}
управления, состоящая в замене оптимального управления на любое допустимое
значение управления лишь на малом интервале времени. Ширина этого интервала
есть параметр вариации, по которому можно дифференцировать. Если на правый
конец траектории никаких ограничений не наложено (т.н. задача {\it со свободным
правым концом}), то целевой функционал на этом однопараметрическом семействе
траекторий должен иметь минимум в нуле, и поэтому его производная по ширине
иголки должна быть неотрицательной. Отсюда легко вытекают условия ПМ.\,
Если же в задаче имеются ограничения на правый конец траектории, то при
таком варьировании полученная траектория, вообще говоря, не будет им удовлетворять.
Чтобы обеспечить попадание конца проварьированной траектории в заданные
ограничения, семейство вариаций должно быть достаточно богатым, и поэтому
приходится рассматривать не одну, а несколько игольчатых вариаций одновременно,
т.н. {\it пакет} игольчатых вариаций, параметром которого является набор
независимых друг от друга ширин иголок.
\ssk

Игольчатые вариации и их пакеты имеют ``родовой недостаток'', состоящий в том,
что ширины иголок по определению неотрицательны, так что получающиеся с их
помощью функции определены лишь на неотрицательном ортанте конечномерного
пространства (точнее, на его пересечении с окрестностью нуля), и для дальнейшего
анализа полученной вспомогательной задачи, соответствующей данному пакету
иголок, требуется проводить ряд специфических построений, выходящих за рамки
стандартного курса анализа и дифференциальных уравнений. Это сразу утяжеляет
доказательство ПМ, делая его практически неподъемным для большинства слушателей.
Поэтому в курсах оптимального управления как правило ограничиваются случаем
задачи со свободным правым концом, где достаточно использовать лишь одну
иголку и никаких дополнительных построений не требуется.
Для изложения доказательства в задачах {\it с ограничениями
на концы траектории}\, желательно, чтобы вспомогательная задача имела вид
гладкой задачи более-менее стандартного (ранее изученного) типа.  \ssk

Одним из приемов, приводящих к гладкой задаче, является т.н. {\it $v-$замена
времени}, предложенная Дубовицким и Милютиным\, в \cite{DM65} (см. также
\cite{Gir, IT, MDO, OPU, develop}), состоящая во введении нового времени
$\t$ и трактовке исходного времени $t=t(\t)$ как фазовой переменной,
подчиненной уравнению
$$ \frac {dt}{d\t}\;=\; v(\t), \qq \mbox{где}\q v(\t)\ge 0
$$
--- еще одно управление. Обратим внимание, что эта замена времени может  быть
не взаимно-однозначной (и это важно!).\,  Ее простейший вариант, когда
зависимость "старого"\, времени от "нового"\, кусочно-линейна (т.е. $v(\t)$
кусочно-постоянна), подробно изложен в \cite{MDO, OPU}. Нетривиальность такой
замены состоит в том, что малые вариации управления $v$ приводят к игольчатым
вариациям (или близким к ним) исходного управления $u,$ и по сути дела $v-$замена
представляет собой такую параметризацию игольчатых вариаций, которая естественным
образом имеет гладкое продолжение на "отрицательные ширины"\,. Прием этот очень
красивый, однако его применение требует хотя и очень прозрачных в идейном плане,
но технически довольно громоздких и не совсем стандартных конструкций, и опыт
преподавания показал, что такое доказательство ПМ все-таки воспринимается
студентами с трудом. \ssk

Еще один способ доказательства ПМ, также предложенный Дубовицким и Милютиным ---
использование вариаций типа {\it скользящих режимов}. (Сами скользящие
режимы были введены Гамкрелидзе в \cite{Gam} для доказательства существования
решения задачи.)\, Этот способ, также красивый, позволяет далекие обобщения
(например, на задачи со смешанными ограничениями, см. \cite{VNIISI, MDO, develop}),
но требует использования весьма нетривиальных (и нестандартных) фактов линейного
и нелинейного функционального анализа, что делает вряд ли оправданным его
применение в классической понтрягинской задаче, тем более в обязательном
курсе лекций.  \ms

Вернемся к игольчатым вариациям как самому простому классу "немалых"\, вариаций.
На протяжении ряда лет на семинаре кафедры Общих проблем управления мехмата МГУ
под руководством В.М. Тихомирова активно обсуждались вопросы преподавания курса
"Вариационное исчисление и оптимальное управление"\,(ВИОУ), в том числе
различные схемы доказательства ПМ. Основное внимание было направлено на
наиболее простое описание конуса концевых вариаций фазовой переменной,
порожденных пакетами игольчатых вариаций управления. Здесь наиболее удачной
нам представляется конструкция, предложенная Магарил-Ильяевым \cite{Mag}.
Доказательство же самого ПМ следует в \cite{Mag} схеме первоначального
доказательства из книги \cite{Pont} (а также \cite{Bolt}).\, Если указанный конус есть все пространство,
то он будет таковым и для некоторого фиксированного пакета иголок, а тогда
соответствующий нелинейный оператор, определенный на неотрицательном ортанте
ширин иголок в данном пакете, обладает т.н. свойством {\it накрывания в
точке} (образ окрестности нуля содержит окрестность образа оптимальной
точки), что приводит к противоречию с оптимальностью исследуемой траектории.\,
Если же указанный конус не совпадает со всем пространством, то применение
конечномерной теоремы об отделимости легко приводит к ПМ.\,
Доказательство накрывания на конусе опирается в \cite{Mag} либо (в случае
кусочно-непрерывного управления) на конусную модификацию  {\it метода Ньютона},
либо (в случае измеримого управления) на довольно тонкий и нетривиальный
топологический факт ---  {\it теорему Брауэра}\, о неподвижной точке. И того,
и другого нет ни в одном стандартном курсе, тем более в прикладных вузах.
(В \cite{Pont, Bolt} оно опиралось на еще более специфическое топологическое
понятие {\it индекса пересечения}\, или  {\it лемму Шпернера}, которые отсутствуют
даже в курсах мехмата.)\, Кроме того, в \cite{Mag}\, рассматривается лишь
случай задачи с фиксированными концами.
\ms

Другой путь доказательства ПМ с помощью игольчатых вариаций, который
представляется нам более предпочтительным, и которому посвящена данная
статья --- это переход к стандартной конечномерной задаче, аргументом которой
являлись бы ширины иголок данного пакета. Идея эта совершенно естественна
и не нова (см. напр. \cite{ATF, Mich, Car} и недавнюю статью \cite{Kor}),
но для ее реализации требуется определить продолжение иголок на отрицательные
ширины. (Как уже говорилось, альтернативой такого продолжения является
указанная выше $v-$замена.)\, В книге \cite[\S 4.2]{ATF}\, было построено
явное продолжение оператора концевого значения фазовой переменной (соответствующего
данному пакету иголок) на полную окрестность начала координат в пространстве
ширин иголок, однако конструкция этого продолжения реализована весьма
громоздко, и поэтому она также вряд ли может быть успешно изложена в общедоступном курсе.
(В работах \cite{Mich, Car} о необходимости указанного продолжения ничего
не говорилось, так что доказательство ПМ в них фактически осталось неполным.)
\ssk

Где-то в середине 2000-х годов на упомянутом семинаре В.М. Тихомирова
авторами данной статьи было предложено использовать возможность продолжения
{\it произвольной}\, гладкой функции, заданной на пересечении неотрицательного
ортанта $\R^k_+$ с окрестностью нуля, на полную окрестность нуля\, с сохранением
гладкости.\, Конечно, здесь можно было бы воспользоваться классической
теоремой Уитни о продолжении гладкой функции с произвольного компакта на
его окрестность (см. \cite{Malg}), но опять-таки, вряд ли было бы оправдано
в данном простом случае обращаться к такой сложной теореме, сама формулировка
которой уже весьма громоздка и требует введения дополнительных понятий.\,
Мы предложили простой способ такого продолжения для произвольной функции,
определенной на $\R^k_+$ и  {\it строго дифференцируемой}\, лишь в нуле.
Этого для наших целей достаточно. (Если функция  {\it непрерывно дифференцируема}\,
всюду на окрестности нуля в $\R^k_+\,,$ что также достаточно для наших целей,
то продолжить ее на полную окрестность нуля можно покоординатным образом,
что сводится к очевидному одномерному случаю. Этот способ использовался
в работе \cite{Kor}, которая является первой известной нам публикацией
с явным построением требуемого продолжения;\, подробнее см. ниже.)\,
\ssk

При наличии подобного продолжения, так или иначе построенного, отпадает
необходимость определять иголки на отрезках отрицательной ширины, вводить
и описывать конус концевых вариаций, а доказательство ПМ состоит в следующем.
Для каждого пакета игольчатых вариаций мы получаем вспомогательную (или
"присоединенную") {\it гладкую конечномерную задачу}\, в пространстве ширин
иголок данного пакета, в которой начало координат есть точка локального
минимума. Применяя стандартное правило множителей Лагранжа, получаем
"конечнозначный"\, ПМ, в котором условие максимума функции Понтрягина
выполнено {\it для данного набора иголок}. Таким образом, мы имеем
{\it семейство конечномерных задач}, соответствующих всевозможным пакетам
иголок, и в каждой задаче свой конечнозначный ПМ.\,
Далее, множество наборов множителей Лагранжа в каждой такой задаче является
конечномерным компактом (это простой факт), а семейство этих компактов
оказывается {\it центрированным} (любое их конечное число имеет непустое пересечение).
Поэтому и все это семейство компактов имеет непустое пересечение, а любой
элемент этого пересечения (т.е. некоторый набор множителей Лагранжа)
обеспечивает "универсальный"\, ПМ, в котором условие максимума выполнено
{\it для любого набора иголок}, т.е. для любых моментов времени и любых
допустимых значений управления.  \ssk

Это доказательство использует лишь стандартные факты матанализа 1 курса,
дифференциальных уравнений, и правило множителей Лагранжа для гладкой
конечномерной задачи с ограничениями равенства и неравенства (которое
предполагается уже прочитанным в курсе ВИОУ к моменту изложения ПМ).\,
Как и в других доказательствах, мы также используем некоторый топологический
факт --- свойство центрированной системы компактов, но это есть
одно из эквивалентных определений компакта, и оно не должно представлять
никаких трудностей для усвоения, тем более в конечномерном случае.
\ssk

По сути дела, наше доказательство следует той же схеме, что и в книге
\cite{ATF};\, отличие лишь в более простом способе продолжения функций на
пакете иголок и несколько более простом изложении дифференциальных свойств
этих функций. От \cite{Mag}\, оно отличается тем, что вместо писания
концевого конуса и доказательства теоремы о конусном накрывании мы переходим
к конечномерной задаче, применяем правило множителей Лагранжа и затем
пользуемся свойством центрированной системы компактов. По нашему мнению,
этот способ доказательства ПМ для классической понтрягинской задачи
(в ее самой общей постановке) в некоторых отношениях проще тех, которые
обычно излагаются в курсах лекций, и поэтому он заслуживает внимания.
\ms

Отметим также, что схема доказательства ПМ, основанная на введении
семейства т.н. {\it присоединенных}\, гладких задач, выписывания условий
стационарности в каждой из них, и затем в использовании центрированной
системы компактов (или, более общо, проективной системы компактов), была
предложена А.Я. Дубовицким и А.А. Милютиным (см. \cite{DM69, Mil70})\,
и успешно применялась не только в задаче понтрягинского типа, но и в более
широком классе задач оптимального управления, включающем как фазовые, так
и смешанные ограничения. В данной статье мы не рассматриваем эти обобщения,
отсылая читателя к работам \cite{VNIISI, MDO, develop}.
\ms

Перейдем к детальному изложению. \\

%\end{document}

%%-------------------------
\section{Общая  задача понтрягинского типа, формулировка принципа максимума}

На отрезке времени $[t_0,t_1],$ априори не фиксированном, будем рассматривать
следующую задачу оптимального управления, которую назовем
{\it общей (или канонической) задачей понтрягинского типа}:
\begin{equation}\label{4}
\dot x(t)=f(t,x(t),u(t)), \qquad u(t)\in U,
\end{equation}
\begin{equation}\label{3}
K(t_0,x(t_0),t_1,x(t_1))=0,
\end{equation}
\begin{equation}\label{2}
F(t_0,x(t_0),t_1,x(t_1))\le 0,
\end{equation}
\begin{equation}\label{1}
J =\; F_0(t_0,x(t_0),t_1,x(t_1))\to \min.
\end{equation}
Здесь $x\in\R^n$ --- фазовая переменная, $u\in\R^r$ --- управление,
концевые функции $K,\, F$ имеют размерности $d(K),\, d(F).$ (Отметим, что
такое обозначение размерностей векторов удобно, так как позволяет
"сэкономить"\,  буквы и освобождает от запоминания размерностей.) \ssk

Предполагается, что $F_0\,,F,\, K\;-$ функции класса $C^1,$ а $f$ непрерывна
вместе с производными $f_t$ и $f_x\,.$ Множество $U \subset \R^r$ произвольно.
Задачу (\ref{4})--({\ref{1}) для краткости назовем {\em задачей} $A$.
\ms

{\bf Замечание 1.}\, Если говорить более точно, свойства функции $f$ надо
предполагать выполненными на множестве ${\cal Q} \times U,$ где ${\cal Q}\;-$
некоторое открытое множество в $\R^{n+1},$ а свойства функций $F_0\,, F,\, K$
на некотором открытом множестве ${\cal P}\subset\R^{2n+2}.$ Как правило, это
всегда подразумевается без явного указания множеств ${\cal Q}$ и ${\cal P}.$
Мы также не будем отвлекаться на эти несущественные детали.\,
\ssk

{\bf Замечание 2.}\, Как отмечалось в \cite{Pont}, множество значений управления
$U$ не обязательно считать лежащим в конечномерном пространстве;\, в общем
случае оно может быть произвольным хаусдорфовым топологическом пространством.
Все нижеследующие рассуждения при этом остаются без изменений. Однако такое
обобщение будет отвлекать читателя (и слушателя) на несущественные моменты,
а для приложений оно вряд ли необходимо. Поэтому, как и в книге \cite{Pont},
мы считаем, что $U \subset \R^r.$
\ms

Решение задачи А\, ищется в классе абсолютно непрерывных функций $x(t)$
и измеримых ограниченных функций $u(t).$ Пару функций $w(t)=(x(t), u(t))$
вместе с отрезком их определения $[t_0, t_1]$ будем называть {\it процессом}\,
задачи. Процесс называется {\em допустимым}, если он удовлетворяет всем
ограничениям задачи. При этом условия (\ref{4}) предполагаются выполненными
почти всюду.\, Как обычно, будем говорить, что допустимый процесс
$\hw(t)=(\hx(t), \hu(t))\mid\, t\in [\hat t_0, \hat t_1],$
доставляет {\it сильный минимум,} если $\ex \e>0,$ такое что $J(w) \ge J(\hw)$
для всех допустимых процессов $w(t)=(x(t), u(t))\mid\, t\in [t_0, t_1],$
удовлетворяющих условиям:
$$
|t_0-\hat t_0|<\varepsilon, \qquad |t_1-\hat t_1|<\varepsilon,\quad
$$ $$
|x(t)-\hat x(t)|<\varepsilon,  %\q |u(t)-\hat u(t)|\,dt <\varepsilon
\qquad  \all t\in [t_0,t_1]\cap [\hat t_0,\hat t_1].
$$ \ssk

Сформулируем теперь {\bf принцип максимума Понтрягина} --- необходимое условие
сильного минимума в задаче А.\, Введем {\em функцию Понтрягина}
$$
H(\p_x,t,x,u)\;=\; \psi_x f(t, x, u),
$$
где $\psi_x$ есть вектор-строка размерности $n$ (зависимость $H$ от  $\psi_x$
иногда будем опускать), и {\em концевую функцию Лагранжа}
$$
l(t_0,x_0,t_1,x_1)=\; (\alpha_0 F_0 + \alpha F + \beta K)(t_0,x_0,t_1,x_1),
$$
где $\alpha_0$ --- число, $\alpha,\, \beta$ --- вектор-строки тех же размерностей,
что и $F,\,K$ соответственно (зависимость $l$ от $\alpha_0\,,\, \alpha,\,\beta$
мы опускаем). \ms

Пусть $w=(x(t), u(t))\mid\,  t\in [t_0, t_1]$ --- допустимый процесс задачи
A.\, Будем говорить, что для него выполнен {\em принцип максимума Понтрягина},\,
если существуют число $\alpha_0\,,$ векторы-строки $\alpha\in\R^{d(F)},$
$\beta\in \R^{d(K)}$ и абсолютно непрерывные функции $\psi_x(t),\; \psi_t(t)$
размерности $n,\,1$ соответственно (где $x$ и $t$ --- индексы, а не
обозначения производных), такие что

\begin{itemize} \item[(i)] $\alpha_0\ge0$, $\alpha\ge0$;

\item[(ii)] $\alpha_0+|\alpha|+|\beta|>0$;

\item[(iii)] $\alpha F(t_0,x(t_0),t_1,x(t_1))=0$;

\item[(iv)] $-\dot\psi_x(t)=H_x(\psi_x(t),t,x(t), u(t)),\qq
-\dot \psi_t(t)=H_t(\psi_x(t),t,x(t), u(t)),$

\item[(v)]

$\psi_x(t_0)= l_{x_0}(t_0,x(t_0),t_1,x(t_1)),\qquad
\psi_x(t_1)= -l_{x_1}(t_0,x(t_0),t_1,x(t_1)); \\[6pt]
\psi_t(t_0)= l_{t_0}(t_0,x(t_0),t_1,x(t_1)),\qquad
\psi_t(t_1)= -l_{t_1}(t_0,x(t_0),t_1,x(t_1)),$

\item[(vi)] $H(\psi_x(t),t,x(t), u(t)) + \psi_t(t) =0$ \q\,
для почти всех $t\in[t_0,t_1],$

\item[(vii)] $H(\psi_x(t),t,x(t), u')+ \psi_t(t) \le 0$\q
для всех\, $t\in[t_0,t_1]$ и всех\, $u'\in U$.

\end{itemize}

Условия $(i)-(v)$  называются условиями неотрицательности, нетривиальности,
дополняющей нежесткости, сопряженными уравнениями, и условиями трансверсальности,
соответственно. \, Условие $(vi)$ не имеет пока стандартного названия;\,
"в рабочем порядке"\, будем называть его {\it законом изменения энергии},
так как из него и сопряженного уравнения для $\p_t$ следует уравнение для
функции $H,$ которая в механических задачах как правило имеет
смысл энергии системы:
$$
\dot H \;=\; H_t\, \qq \mbox{или} \qq \frac{dH}{dt}\;=\; \frac{\pd H}{\pd t}\;.
$$
(В случае, когда управляемая система автономна, т.е. $f = f(x,u)$ не зависит
от $t,$ получаем {\it закон сохранения энергии}: $\dot H =0,$ т.е.
$H = \const).$ \ssk

Из условий $(vi)$ и $(vii)$ вытекает {\it условие максимума}\, функции Понтрягина:
$$
\max_{u' \in U}\, H(\psi_x(t),t,x(t), u') \; =\; H(\psi_x(t),t,x(t), u(t))
\q\, \mbox{для почти всех}\;\; t\in[t_0,t_1],
$$
которое и дало всей совокупности условий $(i)-(vii)$ название {\em принцип
максимума Понтрягина}.

%Заметим, что первое из двух условий $(vi)$ фактически определяет сопряженную
%переменную $\psi_t(t)$ с помощью функций $\psi_x(t)$, $x(t)$ и $u(t)$ (как
%абсолютно непрерывную функцию, которая на множестве полной меры совпадает с
%измеримой функцией $-H(t,x(t), u(t))=-\psi_x(t)f(t,x(t), u(t))$).

Отметим, что обозначения $\psi_x(t)$ и $\psi_t(t)$ для сопряженных переменных
предложены А.Я. Дубовицким и А.А. Милютиным. Их удобство быстро выясняется
при решении конкретных задач с многими фазовыми переменными. \ms

{\bf Замечание 3.}\, Можно показать, что уравнение для функции $\p_t$
(т.е. для $-H)$ вытекает из остальных условий ПМ, т.е. оно не является
независимым. (Ниже мы покажем это для автономной задачи.)\,
Тем не менее, его включение в число условий ПМ вполне оправдано не только
связью с принципами механики, но и тем фактом, что во многих задачах
удобно использовать это условие непосредственно, в "готовом виде"\,,
не занимаясь его выводом из остальных условий ПМ. \ms

Необходимые условия сильного минимума даются следующей теоремой. \ms

\begin{theorem}\label{th1}
Если процесс $\hw=(\hx(t), \hu(t))\mid  t\in [\hat t_0, \hat t_1])$ доставляет
сильный минимум в задаче A, то для него выполнен принцип максимума Понтрягина.
\end{theorem}

Сначала мы докажем эту теорему для случая, когда отрезок времени фиксирован,
а управляемая система автономна. (Последнее, впрочем, не обязательно.)
\bs

%%-------------------------------
\section{Принцип максимума для задачи на фиксированном отрезке времени}
%формулировка и доказательство

Рассмотрим следующую задачу В:
\begin{equation}\label{B4}
\dot x(t)= f(x(t),u(t)), \qquad u(t)\in U, \qq t\in[t_0,t_1],
\end{equation}
\begin{equation}\label{B3}
K(x(t_0),\,x(t_1))=0, \q
%\end{equation}
%\begin{equation}\label{B2}
F(x(t_0),\,x(t_1)) \le 0, \q
\end{equation}
\begin{equation}\label{B1}
J =\; F_0(x(t_0),\,x(t_1))\to \min.
\end{equation}
Здесь отрезок $[t_0,t_1]$ фиксирован, функции $F_0$, $F$ и $K$  непрерывно
дифференцируемы, а функция $f$ непрерывна вместе с производной  $f_x$.\,
Минимум ищется среди всех пар (процессов)
$w=(x,u) \in W: = AC^n[t_0,t_1]\times L_\infty^r[t_0,t_1],$
где $AC^n[t_0,t_1]$ есть пространство абсолютно непрерывных функций
размерности $n,$ а $L_\infty^r[t_0,t_1]\;- $ пространство измеримых
ограниченных функций размерности $r.$
\ssk

По определению, допустимый процесс $\hw=(\hx,\hu)$ доставляет {\em сильный
минимум}, если существует $\varepsilon>0$ такое, что $J(w)\ge J(\hat{w})$ для
всех допустимых пар $w=(x, u),$ удовлетворяющих  условию $||x-\hat x||_C <\e.$
Другими словами, сильный минимум в задаче В\, это локальный минимум в
пространстве $W$ относительно полунормы $||w||' = ||x||_C\,.$ \ssk

Как и выше, вводим функцию Понтрягина $H(\p,x,u)\,=\, \psi f(x, u),$
где $\psi$ есть $n-$мерная вектор-строка (в общей задаче А\, она обозначалась
$\psi_x),$ и  концевую функцию Лагранжа
$l(x_0,\,x_1)=\; (\alpha_0 F_0 + \alpha F + \beta K)(x_0,\,x_1).$
\ssk

Условия ПМ для допустимого процесса $w=(x, u)$ задачи В\, состоят в следующем:\,
существуют число $\alpha_0\,,$ векторы-строки $\alpha\in\R^{d(F)},$
$\beta\in \R^{d(K)},$ абсолютно непрерывная функция $\psi(t)$ размерности
$n\,$ и константа $c\in\R,$ такие что

\begin{itemize} \item[(i')] $\alpha_0\ge 0$, $\alpha\ge 0$;

\item[(ii')] $\alpha_0+|\alpha|+|\beta|>0$;

\item[(iii')] $\alpha F(x(t_0),\,x(t_1))=0$;

\item[(iv')] $-\dot\psi(t)=\; H_x(\p(t), x(t), u(t)),$

\item[(v')]

$\psi(t_0)= l_{x_0}(x(t_0),\,x(t_1)),\qquad
\psi(t_1)= -l_{x_1}(x(t_0),\,x(t_1)),$

\item[(vi')] $H(\p(t), x(t), u(t)) =\; c$ \q\, почти всюду на $[t_0,t_1],$

\item[(vii')] $H(\p(t), x(t), u')\lee c$ для всех\,
$t\in[t_0,t_1]$ и всех\, $u'\in U.$

\end{itemize}

Заметим, что здесь вместо сопряженного уравнения на $\psi_t$ мы пишем
условие постоянства функции Понтрягина вдоль оптимального процесса
({\it "закон сохранения энергии"\,},\, как и полагается в консервативной
системе).
\ms

Имеет место следующая

\begin{theorem}\label{th2}\,\,
Если процесс $\hat w=(\hat x, \hat u)$ доставляет сильный минимум в задаче
$B,$ то для него выполнен принцип максимума $(i')-(vii').$
\end{theorem}

Доказательство\, мы проведем для случая, когда управление $\hat u(t)$
кусочно непрерывно. Общий случай, когда оптимальное управление измеримо,
технически более сложен и, как правило, не рассматривается в стандартных
курсах.
%(Ниже мы укажем, что изменится в общем случае измеримой ограниченной $\hu(t))$.
\ssk

Предварительно изучим некоторые дифференциальные свойства управляемой системы
\beq{xua1}
\dx = f(x,u)\qq % x(t_0)= x_0\,
\eeq
на произвольном ее решении $(x(t),\,u(t))$ на  отрезке $[t_0,t_1].$
Будем рассматривать управления, равномерно ограниченные некоторой константой:
$|u(t)| \le M.$ Из теории ОДУ известно (впрочем, это несложно доказать и
непосредственно), что если при этом $||u-\hu||_1$ мала, и начальное значение
$x_0=x(t_0)$ мало отличается от $\hx(t_0),$ то решение $x(t)$ уравнения (\ref{xua1})
всегда существует на всем отрезке $[t_0,t_1]$ и мало отличается от $\hx(t)$
в норме пространства $C[t_0,t_1].$  \ssk

Введем для процесса $(x(t),\,u(t))$ с кусочно непрерывным управлением
\bs

%%------------------------------------------------------------
{\bf Понятие элементарной игольчатой вариации}.\, Будем считать, что
управление $u(t)$ непрерывно слева во всех точках полуинтервала $(t_0,t_1]$
и непрерывно справа при $t=t_0\,.\,$
Фиксируем любую точку $\theta\in (t_0,t_1)$ и любое значение $v\in U.$
Для каждого достаточно малого  $\e\ge 0$ определим на отрезке $[t_0,t_1]$
новое управление
$$
u_\e(t)\,= \;
\left\{\begin{array}{l}
\; v,\quad\;\; \mbox{если}\quad t\in (\theta -\e,\,\theta],\\[4pt]
 u(t),\;\;\; \mbox{если}\quad t\notin (\theta -\e,\,\theta].
\end{array}\right.
$$
(При $\e=0$ полуинтервал $(\theta -\e,\,\theta]$ пуст, и управление не меняется.)\,
Отметим, что это управление по-прежнему непрерывно слева.\, Семейство
функций $u_\e(t)$ при $\e\to 0+$ называется {\it элементарной игольчатой
вариацией}\, исходного управления $u(t)$ с параметрами $(\theta,v).$
\ssk

Пусть вектор $a\in\R^n$ мало отличается от $ x_0= x(t_0)\,.$ Определим
на отрезке $[t_0,t_1]$ фазовую переменную $x_\e(t)$ как решение задачи
Коши для управляемой системы (\ref{xua1}):
\beq{xea1}
\dx_\e =\, f(x_\e,\, u_\e), \qq x_\e(t_0) =a.
\eeq
Если $\e\ge 0$ достаточно мало и $a$ достаточно близко к {$ x_0\,,$} то,
как уже было сказано, решение задачи (\ref{xea1}) существует на всем отрезке
$[t_0,t_1],$ единственно, и непрерывно в норме пространства $C[t_0,t_1]$
зависит от пары $(a,\e).$ Примем это как установленный факт. Значение
фазовой переменной $x_\e(t)$ в точке $t_1$ обозначим через $P(a,\e).$
\ssk

Таким образом, имеем отображение  \vad
$$
P:\; (a,\e)\; \longmapsto\; x_\e(t_1) \in \R^n,
$$
которое определено и непрерывно на $\W(x_0) \times [0,\e_0),$
где $\W(x_0)$ есть окрестность точки $x_0 =x(t_0) \in \R^n,\,$ а $\e_0 >0.$
При этом очевидно $P( x_0,0)=   x(t_1)=:  x_1$.
\ms

Установим, что отображение $P$ в точке $(a,\e)= (x_0,\,0)$ имеет производную
по $a$ и правую производную по $\e,$ которые далее будем обозначать как
$P_a(x_0,0)$ и $P_\varepsilon^+(x_0,0)$  соответственно. \ssk

Для нахождения производной по $a$ в некотором направлении $\bar a$
надо, как известно, решить {\it уравнение в вариациях}\,
\beq{urvar1}
\dot{\bx} =\, f_x(x,\, u)\,\bx, \qq \bx(t_0)= \bar a,
\eeq
При этом $P_a(x_0,0)\,\bar a\, =\, \bx(t_1).$ Обратим внимание, что если
$x(t)$ мы будем мало менять относительно нормы пространства $C[t_0,t_1],$
а $u(t)$ мало менять относительно нормы пространства $L_1[t_0,t_1]$
(сохраняя его равномерную ограниченность априорной константой), то
матрица $f_x(x, u)$ будет мало меняться относительно нормы пространства
$L_1[t_0,t_1],$ а тогда, как известно, решение уравнения (\ref{urvar1})
будет мало меняться относительно нормы $C[t_0,t_1].$
%(См. лемму \ref{lem8} в Приложении.)
\ms

Для нахождения правой производной $P_{\e}^+(x_0,0)$ при фиксированном $x_0$
будем рассуждать так. Каждому $\e> 0$ соответствует свое решение $x_\e(t)$
уравнения (\ref{xea1});\,\, при $\e=0$ имеем невозмущенное решение $x(t).$
Если $t<\theta,$ то при достаточно малых $\e$ оба решения в этой точке совпадают.
Оценим разность $\D x_\e = x_\e - x$ в точке $\theta.$ Имеем  \vad
$$
\begin{array}{l}
x(\theta)\, =\; x(\theta-\e)\, +\, \e\,f(x(\theta), u(\theta)) + o(\e),  \\[6pt]
x_\e(\theta) =\; x(\theta-\e)\, +\, \e\,f(x(\theta), v)\; +\; o(\e),
\end{array}
$$
откуда $\D x_\e(\theta) =\, \e\cdot \D f(\theta,v) + o(\e),\;$ где
$\D f(\theta,v) = f(x(\theta), v) - f(x(\theta), u(\theta)).$
\ms

Положим $\dis \bx(t) = \lim_{\e\to 0+}\, \frac{\D x_\e(t)}{\e}\,.\;$
Тогда $x_\e(t) = x(t) + \e\,\bx(t) + o(\e).$ \\[8pt]
При этом левее точки $\theta$ имеем очевидно $\bx(t) =0,$ в самой этой
точке $\bx(\theta) = \D f(\theta,v),$ а на отрезке $[\theta,\,t_1]$ обе
функции $x$ и $x_\e$ удовлетворяют уравнению (\ref{xua1}) с одним и тем
же невозмущенным управлением, но с разными начальными значениями в точке
$\theta.$ Поэтому, как и раньше, $\bx$ удовлетворяет на отрезке $[\theta,\,t_1]$
уравнению в вариациях  \vad
\beq{urbx1}
\dot{\bx} =\, f_x(x,\, u)\,\bx \qq
\mbox{с начальным значением}\q
\bx(\theta)= \D f(\theta,v),
\eeq
и тогда $P_{\e}^+(x_0,0) = \bx(t_1).$ \ms

Отметим, что $P_\e^+(x_0,\,0)$ зависит не только от параметров иголки $(\theta,v),$
но и от базовой пары $(x(t), u(t)),$ удовлетворяющей уравнению (\ref{xua1}).
Если последняя мало меняется в норме $||x||_C + ||u||_1$ (при равномерно
ограниченных управлениях) и пара значений $(\theta, u(\theta))$ также
меняется мало, то матрица $f_x(x(t), u(t))$ мало изменится в норме $L_1\,,$
начальный момент $\theta$ и начальное значение $\bx(\theta)$ также
изменятся мало, и тогда решение задачи Коши (\ref{urbx1}) мало изменится
в норме $C,$ а значит, мало изменится и его концевое значение $x(t_1).$
\ms

%%-----------------------------------------------------
Итак, доказана

\ble{lem1}
В точке $(a,\e)= (x_0,\,0)$ отображение $P$ имеет производную по $a$
и правую производную по $\e,$ которые  выражаются формулами:
$$
\all \bar a \in \R^n \qq P_a(x_0,\,0)\,\bar a\, =\, \bx(t_1),
$$
где $\bx(t)$ есть решение задачи Коши (\ref{urvar1}) на отрезке $[t_0,\,t_1],$
а правая производная $\, P_{\e}^+(x_0,0) = \bx(t_1),\,$ где $\bx(t)$ есть
решение задачи Коши (\ref{urbx1}) на отрезке $[\theta,\,t_1].\;$ \ssk

Обе эти производные непрерывно зависят от изменения пары $(x, u)$ в норме
пространства $C \times L_1$ (при условии равномерной ограниченности
управлений) и пары $(\theta, u(\theta))\in \R^{1+r}.$
\ele
\ssk

Если теперь в качестве базового управления взять функцию $u_\e(t)$ при малом
$\e>0,$ в качестве начального значения взять вектор $a,$ близкий к $x(t_0),$
а иголку строить в точке $\theta -\e,$ то разность $u_\e -u$ будет мала
в норме $L_1[t_0, t_1],$ разность $x_\e -x$ мала в норме $C[t_0, t_1],$
и в силу непрерывности слева функции $u(t)$ новое базовое значение
$u_\e(\theta -\e) = u(\theta -\e)$ будет близко к старому $u(\theta).$
Тогда и решение задачи Коши (\ref{urbx1}) на отрезке $[\theta -\e,\,t_1],$
рассматриваемое на отрезке $[\theta,\,t_1],$ будет равномерно близко к
старому решению.\, Поэтому из леммы \ref{lem1} вытекает
\ms

{\bf Следствие.}\, {\it Оператор $P,$ соответствующий элементарной
игольчатой вариации в точке $\theta,$ при любых значениях $a,$ достаточно
близких к $x(t_0),$ и любых достаточно малых $\e \ge 0$ имеет в точке
$(a,\e)$ производную $P_a(a,\,\e)$ и правую производную $P_\e^+(a,\,\e),$
которые непрерывно зависят от пары $(a,\e).$ }
\ms

Конкретные выражения этих производных нам потребуются лишь для точки
$(a,\e) = (x_0,0);\,$ они приведены выше.
\bs

{\small
{\bf Замечание 4.}\, В случае, когда базовое управление измеримо, это
следствие уже не имеет места даже если $\theta$ есть {\it точка Лебега}\,
функции $u(t).$ Дело в том, что в общем случае начальное значение для
уравнения (\ref{urbx1}) дается формулой
$$
\bx(\theta) =\; f(x(\theta), v) -\; \lim_{\d \to 0}\,
\frac 1\d \int_{\theta -\d}^\theta f(x(\theta), u(\t))\,d\t,
$$
если указанный предел существует. В самой точке $\theta$ он существует и
равен $f(x(\theta), u(\theta)),$ так что по-прежнему $\bx(\theta)= \D f(\theta,v),$
но в сдвинутой точке $\theta -\e\,$ при $\e>0$ этот предел может сильно измениться,
и поэтому непрерывности правой производной по $\e$ мы не получим.\,
Можно гарантировать существование указанных производных лишь в данной 
точке $(x_0,0)$ и непрерывность самого отображения $P$ в ее окрестности.
При этом мы выходим за рамки гладких задач и приходится использовать 
тонкие топологические факты типа теоремы Брауэра. В этом смысле возможности
стандартных игольчатых вариаций довольно ограничены.\, Ситуацию все же
можно, по-видимому, исправить, если в качестве $\theta$ рассматривать
{\it точки аппроксимативной непрерывности}\, функции $u(t)$ и полагать
$u_\e(t) =v$ не на интервале $(\theta -\e,\, \theta),$ а на его пересечении
с тем множеством, вдоль которого $u(t)$ непрерывна в точке $\theta.$
Но и тогда мы получим лишь строгую дифференцируемость при $\e=0,$ но
не существование и непрерывность производной при $\e> 0.$ Кроме того,
такой способ потребует довольно сложных конструкций и глубоких
фактов теории меры, что никак не годится для лекционного курса.

Тем не менее, даже в случае измеримой функции $u(t)$ гладкость отображения
$P$ можно получить относительно просто, если игольчатую вариацию делать не
путем {\it замены}\, данного управления $u(t)$ на выбранное значение $v$
на интервале времени длины $\e$ около точки $\theta,$ а путем "расшивки"\,
этой точки до отрезка длины $\e$ и {\it вставки}\, управления $v$ на этом
отрезке.
(При этом весь отрезок $[t_0, t_1]$ увеличится до отрезка $[t_0,\, t_1 +\e].)$
По сути, это и будет результат упомянутой выше $v-$замены времени.
(Букву $v$ мы здесь вынужденно употребили в двух разных смыслах.) }
\bs

Далее нам потребуется дифференцировать скалярную функцию вида
$G(a,\e) = g(a,\, P(a,\e)),$ где $g(x_0,\, x_1)\;-$ дифференцируемая функция в окрестности
точки $(\hx_0,\, \hx_1)$ пространства $\R^{2n}.$ Это удобно рассмотреть как
отдельное свойство.\,
\ms

Пусть опять пара $(x(t), u(t))$ есть решение уравнения (\ref{xua1}).\,
Введем для нее липшицеву функцию $\p(t)$ (вектор-строку размерности $n)$
как решение задачи Коши
\beq{adeq}
\dot\p =\; -\p\,f_x(x,u), \qq \p(t_1) =\, -g_{x_1}(x_0,\,x_1).
\eeq
Отметим, что $\p(t)$ не зависит от игольчатой вариации и ее параметров $(\theta, v).$
\ms

%%------------ Лемма 2 --------------------------------------
\ble{lem2}\q
Производная функции $G$ по переменной $a$ и правая производная этой функции
по $\varepsilon$ в точке $(x_0,\,0)$ выражаются формулами:
$$
\phantom{tttttttttt}  G_a(x_0,\,0)\,\bar a\, =
\,(g_{x_0}(x_0,x_1)- \p(t_0))\,\bar a,
\qq \all \bar a \in \R^n,
$$ $$
G_\e^+(x_0,\,0) =\; -\p(\theta)\cdot \D f(\theta, v).
$$
\ele

\Proof Пусть $\bx$ есть любое решение уравнения в вариациях
$\dot{\bx} =\, f_x(x,\, u)\,\bx\,$ на некотором отрезке $[t', t''],$
неважно с каким начальным значением. Нетрудно видеть, что произведение
$\p(t)\,\bx(t)$ постоянно на этом отрезке. Действительно, его производная
$$
\frac{d}{dt}(\psi\, \bx)=\; \dot\psi\,\bx\,+ \,\psi\dot{\bx}\; =\;
-\,\psi f_x\,\bx\, +\, \psi f_x\,\bx =0.
$$
(Этот факт верен для решения любого линейного уравнения $\dot{\bx} = A(t)\,\bx$
и {\it сопряженного}\, к нему $\dot \psi = -\psi\, A(t)$ с любой интегрируемой
матрицей $A$.) \ctd
\ms

Возьмем теперь произвольное $\bar a \in \R^n,$ и пусть $\bx$ есть решение
уравнения (\ref{urvar1}). Тогда, с учетом леммы \ref{lem1},
$$
G_a(\hx_0, 0) = \; g_{x_0}\,\bar a\,+\, g_{x_1}\,P_a(\hx_0,0)\,\bar a\; =
\; g_{x_0}\,\bar a\,+\,  g_{x_1}\,\bx(t_1) \;=
$$ $$
\, = g_{x_0}\,\bar a\, -\,\p(t_1)\,\bx(t_1) \;=\;
g_{x_0}\,\bar a\, -\,\p(t_0)\,\bx(t_0) \;=\;(g_{x_0}- \p(t_0))\,\bar a.
$$
Аналогично:
$$
G_\e^+(\hx_0, 0) = \; g_{x_1}\,P_\e^+(\hx_0,0)\;=\; -\p(t_1)\,\bx(t_1)\;=
$$ $$
=\; -\p(\theta)\,\bx(\theta) \;=\; -\p(\theta)\cdot \D f(\theta, v).
\qq\qq\qq \mbox{\ctd}
$$ \ssk

Итак, мы изучили дифференциальные свойства элементарной игольчатой вариации.
Однако, как известно, для доказательства ПМ в задаче (\ref{4})--(\ref{1})
одной игольчатой вариации недостаточно. Поэтому определим
\ms

%%--------------------------------------------------
{\bf Понятие пакета игольчатых вариаций.}\,\, Рассмотрим опять базовую пару,
которую обозначим теперь $(\hx(t), \hu(t)),$ с кусочно непрерывным управлением
$\hu(t),$ которое считаем непрерывным слева во всех точках полуинтервала
$(t_0,t_1]$ и непрерывным справа при $t=t_0\,.\,$ \ssk

Пусть $\cal N$ есть произвольный конечный набор пар $(\theta_i,\,v_i),$
$i=1,\ldots,s,$ где $\theta_1\le\ldots\le\theta_s\;-$ любые точки из интервала
$(t_0,t_1),$ а $v_1\,, \ldots,\, v_s\;-$ любые точки из множества $U.$
Пусть задан вектор $\varepsilon=(\varepsilon_1,\ldots,\varepsilon_s)$
с малыми неотрицательными компонентами $\varepsilon_i\,.$ Построим
вариацию управления, которая принимает значения $v_i$ на полуинтервалах
$\Delta_i$ длины $\e_i$ около точек $\theta_i\,.$
\ssk

Определим эти полуинтервалы следующим образом. Если все точки $\theta_i$
различны, то полагаем $\D_i = (\theta_i -\e_i,\, \theta_i],$ $i=1,\ldots,s.$
Если же среди $\theta_i$ есть повторяющиеся, то для каждого набора повторяющихся
значений $\theta_k = \ld = \theta_{k+p}$ откладываем последовательно примыкающие
полуинтервалы длины $\e_{k+i}$ влево от точки $\theta_k\,,$ т.е. полагаем
$$
\D_k = (\theta_k -\e_k,\, \theta_k], \q\;
\D_{k+1} = (\theta_k -\e_k -\e_{k+1},\; \theta_k- \e_k], \q\; \mbox {и т.д.},
$$
а для всех не повторяющихся значений $\theta_i$ полагаем, как и прежде,
$\D_i = (\theta_i -\e_i,\, \theta_i].$ Итак, для данного набора $\cal N$ и
вектора $\e\ge 0$ мы определили полуинтервалы $\D_i$ длины $\e_i\,.$
Если вектор $\varepsilon$ достаточно мал, эти полуинтервалы не пересекаются
и все находятся на отрезке $[t_0,t_1].$ Обозначим через $\t_i(\e)$ правый
конец полуинтервала $\D_i\,.$ Ясно, что он непрерывно зависит от вектора
$\e.$  \ssk

Определим теперь управление
$$
u_\e(t)\,=\;
\left\{\begin{array}{l}
v_i\,,\quad \; \mbox{ если}\quad t\in \Delta_i\,,\quad i=1,\ldots,s, \\[6pt]
\hu(t),\quad \mbox{если}\quad t \in [t_0, t_1] \setminus \bigcup_{i=1}^s \Delta_i\,.
\end{array}\right.
$$
Семейство функций $u_\e(t)$ при $\e\to 0+$ называется {\it пакетом
игольчатых вариаций}\, исходного управления $\hu(t).$ \ssk

Подставив $u_\e(t)$ в управляемую систему (\ref{B4}), получим задачу Коши
(\ref{xea1})
%\beq{xea2} \dx_\e =\, f(x_\e,\, u_\e), \qq x_\e(t_0) =a, \eeq
с некоторым начальным значением $a.$ Как и в случае элементарной игольчатой
вариации, из теории ОДУ следует, что при достаточно малых $\e\ge 0$ и при
начальных значениях $a,$ достаточно близких к $\hx_0 =\hx(t_0),$ решение
уравнения (\ref{xea1}) существует, единственно, и непрерывно зависит от
набора $(a,\e).$
%Примем это как установленный факт.
\ssk

Таким образом, имеем отображение
$$
P:\; (a,\e)\; \longmapsto\; x_\e(t_1) \in \R^n,
$$
которое определено и непрерывно на $\W(\hx_0) \times (\calO \cap \R^s_+),$
где $\W(\hx_0)$ есть окрестность точки $\hx_0 \in \R^n,$ через $\calO$ обозначена
окрестность нуля в $\R^s,$ а $x_\e(t_1)$ есть значение решения (\ref{xea1})
в конечной точке $t=t_1\,.$

Из леммы \ref{lem1} следует, что отображение $P$ не только непрерывно,
но и является гладким в следующем смысле.

%%----------- Лемма 3 --------------------
\ble{lem3}
В любой точке $(a,\e)$ своей области определения отображение $P$ имеет
производную по $a$ и правую производную по каждому $\e_i\,,$ которые
непрерывно зависят от пары $(a,\e).$
\ele

\Proof Существование указанных производных (и их выражение через конкретные
формулы) установлено в лемме \ref{lem1}. Далее, при изменении пары $(a,\e)$
управление $u_\e$ непрерывно меняется в норме пространства $L_1[t_0, t_1],$
соответствующее ему $x_\e$ непрерывно меняется в норме пространства
$C[t_0, t_1],$ и для каждой $i-$ой иголки пара её "базовых"\, значений
$(\t_i(\e),\, u_\e(\t_i(\e))$ меняется непрерывно при изменении $\e$
в силу непрерывности слева функции $u(t).$ Поэтому в силу леммы 1 производные
отображения $P$ будут непрерывно зависеть от пары $(a,\e)$ в окрестности
точки $(\hx_0\,,0).$ \ssk

Поясним сказанное на примере двух иголок в точке $\theta$ со значениями
$v_1, v_2$ на полуинтервалах $\D_1 = (\theta -\e_1,\, \theta], \,$
$\D_2 = (\theta -\e_1 -\e_2,\, \theta -\e_1].$ По лемме~\ref{lem1} здесь
$P_{e_1}^+$ выражается формулой для обычной элементарной иголки в точке
$\theta,$ а $P_{e_2}^+(\e_1, \e_2) = \bx(t_1),$ где $\bx(t)$ есть решение
задачи Коши (\ref{urbx1}) на отрезке \linebreak
$[\theta -\e_1 -\e_2,\,t_1]\,$ с начальным значением \\[4pt]
\hspace*{10mm}
$ \bx(\theta -\e_1 -\e_2) =\; f(\hx(\theta -\e_1 -\e_2),\, v_2)
- f(\hx(\theta -\e_1 -\e_2),\, \hu(\theta -\e_1 -\e_2)).\,$
\\[4pt]
Обе эти производные непрерывно зависят от $(\e_1,\e_2)$ вплоть до границы
квадранта $\R^2_+\,.$ В частности, $P_{e_2}^+(\e_1, 0)$ соответствует начальному
значению $\bx(\theta -\e_1) = f(\hx(\theta -\e_1),\, v_2) -
f(\hx(\theta -\e_1),\, \hu(\theta -\e_1)),$ а $P_{e_1}^+(0, \e_2)\;- $
начальному значению $\bx(\theta -\e_2) =
f(\hx(\theta -\e_2),\, v_2) - f(\hx(\theta -\e_2),\, \hu(\theta -\e_2)).$
\ctd
\bs

Покажем теперь, что если $\e \in \inter\,\R^s_+$ (кратко пишем $\e>0$), то
существуют обычные частные производные $P_{\e_i}(a,\e),$ совпадающие с
правыми производными. Это одномерный факт, он вытекает из следующего
утверждения, которое является несложным упражнением по матанализу.

\ble{lem4}
Пусть функция $\f: \R_+ \to \R$ при всех $x\ge 0$ имеет правую производную
$\f'_{\mbox{\rm \small пр}}(x),$ которая непрерывна на $\R_+\,.$ Тогда
$\f$ при всех $x>0$ имеет обычную производную, которая совпадает с правой:
$\; \f'(x) = \f'_{\mbox{\rm \small пр}}(x).$
\ele

Доказательство дано в Приложении. \bs

Из этой леммы и леммы \ref{lem3} следует, что отображение $P$ непрерывно дифференцируемо в
области $\e>0,$ точнее на $\W(\hx_0) \times (\calO \cap \inter\R^s_+)).$
Кроме того, его производная имеет предел в тех точках $(a,\e),$ где
$\e \in \pd\R^s_+\,;\,$
%$(a,\e) \in \W(\hx_0) \times (\calO \cap \pd\R^s_+)),$ который
его мы также будем обозначать $P'(a,\e).$ Отсюда, в свою очередь следует,
что $P$ {\it строго дифференцируемо}\, в точке $(\hat x_0,0),$ т.е.
разность $P(a,\e) - P'(\hat x_0,0)(a,\e)$ липшицева в окрестности точки
$(\hat x_0,0)$ (пересеченной, естественно, с $\R^n \times \R^s_+)$ с константой,
стремящейся к нулю вместе с радиусом окрестности. (Понятие строгой
дифференцируемости мы считаем уже изложенным ранее в данном курсе
оптимального управления.) \ms

Таким образом, гладкость отображения $P,$ т.е. непрерывность его частных
производных, сводится к непрерывности производных одной элементарной
иголки в зависимости от базового процесса. Эту идею мы заимствовали
у Магарил--Ильяева \cite{Mag}.  \ms

Ключевой пункт нашего доказательства состоит в том, что отображение $P$ можно
продолжить на "полную"\, окрестность точки $(\hat x_0,0),$ т.е. на произвольные
$a \in \W(\hx_0)$ и произвольные малые значения вектора $\e,$ сохраняя его
строгую дифференцируемость в точке $(\hat x_0,0).$ Это вытекает из следующего
общего утверждения. \ssk

%%----------------------------------------------------------------------
\ble{lem5}
Пусть отображение $P: \R^n \times K \to \R^m,$ где $K$ есть замкнутый
выпуклый конус в $\R^s$ с непустой внутренностью, строго дифференцируемо
в точке $(0,0).$ Тогда его можно продолжить до отображения
$\wt P: \R^n \times \R^s \to \R^m,$ так же строго дифференцируемого
в точке $(0,0).\,$
\ele

Его доказательство мы отнесли в Приложение. (У нас в качестве конуса $K$
выступает неотрицательный ортант $\R^s_+$.)\, Указанное продолжение будем
обозначать той же буквой $P.$ Отметим, что для $K = \R^s_+$ и гладкого
$P$ (как в нашем случае) можно обойтись и без понятия строгой
дифференцируемости, используя лишь стандартное понятие непрерывной
дифференцируемости (см. \cite{Kor}):
\ssk

\ble{lem5a}
Пусть отображение $P: \R^n \times \R^s_+ \to \R^m$ непрерывно дифференцируемо
внутри своей области определения, и его производная непрерывна вплоть до
границы.\, Тогда его можно продолжить до отображения
$\wt P: \R^n \times \R^s \to \R^m,$ непрерывно дифференцируемого на всем
пространстве.
\ele

Доказательство мы также отнесли в Приложение. \ssk

Обратим еще раз внимание, что все предыдущие рассуждения относятся только
к дифференциальным свойствам управляемой системы на пакете игольчатых вариаций\,
и не связаны непосредственно с оптимизационной задачей. Поэтому они и сами
по себе представляют интерес и могут быть использованы также и в других
вопросах теории управления.
\ms

Теперь мы готовы к доказательству ПМ для задачи В.
\bs

%%-------------------------------
\subsection*{Принцип максимума для пакета игольчатых вариаций.}

Пусть процесс $\hat w=(\hat x, \hat u)$ доставляет сильный минимум
в задаче B.\, Зафиксируем некоторый пакет ${\cal N}$ игольчатых вариаций.
С помощью леммы \ref{lem5} (или леммы \ref{lem5a}) продолжим отображение
$P$ на полную окрестность точки $(\hx_0,0),$ где $\hx_0 = \hx(t_0),$
с сохранением строгой дифференцируемости в этой точке (или, соответственно,
непрерывной дифференцируемости в ее окрестности). \ssk

В конечномерном пространстве $\R^n\times \R^s$ с элементами $(a,\e)$
рассмотрим следующую {\em задачу} $Z(\calN)$:
$$
J(a,P(a,\varepsilon)) \to \min,  \quad F(a,P(a,\varepsilon))\le 0,
\quad K(a,P(a,\varepsilon))=0, \quad -\varepsilon\le 0.
$$
Эта задача является сужением задачи В, и поэтому из того, что пара
$(\hat x,\hat u)$ доставляет сильный минимум в задаче B\, вытекает, что
пара $(\hat a= \hx_0,\,\,\hat\e=0)$ доставляет локальный минимум в задаче
$Z(\calN),$ и следовательно, удовлетворяет условиям стационарности.

Выпишем эти условия для точки $(\hx_0,0)$ в задаче $Z(\calN).$
Они состоят в том, что существуют множители Лагранжа
$$
\alpha_0\in\R,\quad \alpha\in {\R^{d(F)}},\quad \beta\in {\R^{d(K)}},
\quad \g\in\R^s,
$$
такие, что выполнены условия
\begin{eqnarray} && \label{stat1a}
\alpha_0\ge 0,\quad \alpha \ge 0, \quad \g \ge 0,\\[4pt]
&& \label{stat2a} \alpha_0+|\alpha|+ |\beta|+ |\g| >0,\quad
\alpha F(\hat x_0,\hat x_1)=0,\\[4pt]
&& \label{stat3a}{\cal L}_{a}(\hat x_0,0)=0,\qquad
{\cal L}_\varepsilon(\hat x_0,0)=0,
\end{eqnarray}
где  \vad
$$
{\cal L}(a,\e)= \big(\alpha_0 J + \alpha F +\beta K\big)(a,P(a,\e))\, -\,\g\,\e\;
= \;  { l} (a,P(a,\e))\,- \g\,\e
$$
-- функция Лагранжа задачи $Z(\calN).$ \ssk

Отметим еще раз, что в отличие от ряда других доказательств ПМ, в т.ч. от
приведенных в \cite{Pont, Bolt, Mag}, мы не строим конус вариаций $\bx(t_1),$
а просто рассматриваем конечномерную задачу, соответствующую данному
пакету игольчатых вариаций, и пользуемся уже известными (и изложенными
ранее в курсе) условиями стационарности для этой задачи. \ms

Расшифруем условия (\ref{stat3a}).\, Пусть липшицева функция $\p(t)$ есть решение
{\it сопряженного уравнения}\, (к уравнению (\ref{urvar1}) вдоль оптимального
процесса $(\hx,\hu)$)
\begin{equation}\label{adj}
\dot \psi =\, -\psi f_x(\hat x(t),\hat u(t))
\end{equation}
с граничным условием  $\;\psi(t_1)= -l_{x_1}(\hat x_0,\hat x_1)\,.$
Обратим внимание, что $\p$ зависит только от набора множителей $(\a_0,\a,\b)$
(поскольку концевая функция $l$ выражается через них) и не зависит
непосредственно от выбранного пакета вариаций. Более того, $\p$ единственным
образом выражается через $(\a_0,\a,\b).$
\ssk

По лемме \ref{lem2}, примененной к функции $g= \calL,$ первое равенство
в (\ref{stat3a}) означает, что $\all \bar a \in \R^n$
$$
{\cal L}_{a}(\hat x_0,0)\,\bar a\, = \;
(l_{x_0}(\hx_0, \hx_1)\, -\,\p(t_0))\,{\bar a}\;=\,0,
$$
откуда следует, что $\p(t_0) = l_{x_0}(\hx_0, \hx_1).$ Таким образом,
функция $\p$ удовлетворяет следующим граничным условиям на концах отрезка:
\beq{tr1}
\p(t_0) = l_{x_0}\,, \qq \p(t_1) = -l_{x_1}\,.
\eeq
Они называются {\it условиями трансверсальности}.
\ms

Второе равенство в (\ref{stat3a}) означает, что для каждого $i$
\beq{lei}
{\cal L}_{\varepsilon_i}(\hx_0,0)= \;
-\p(\theta_i)\cdot \D f(\theta_i, v_i) - \g_i\; = 0,
\eeq
т.е. $\q \p(\theta_i)\cdot \D f(\theta_i, v_i)\, =\;
\p(\theta_i)\Bigl( f(\hx(\theta_i),v_i)) -
f(\hx(\theta_i), \hu(\theta_i)) \Bigr)\, =\, - \g_i \lee 0.$
\bs

Введем функцию Понтрягина $\; H(\p,x,u) = \p\,f(x,u).$ Тогда последнее
условие означает, что выполнено {\it конечнозначное условие максимума}
\beq{usmax}
H(\p(\theta_i),\hx(\theta_i),v_i) \lee H(\p(\theta_i),\hx(\theta_i), \hu(\theta_i)),
\qq  i=1, \ld, s,
\eeq
а сопряженное уравнение (\ref{adj}) может быть записано в виде
\beq{adjh}
\dot\p(t) =\, - H_x (\p(t),\hx(t), \hu(t)).
\eeq
(Условие (\ref{usmax}) можно трактовать так: \, при любом $\theta_i$ из
данного пакета функция $H(\p(\theta_i),\hx(\theta_i),u)$ принимает максимальное
значение по всем значениям управления $u,$ имеющимся в этом пакете в
точке $\theta_i\,,$ при оптимальном $u = \hu(\theta_i).)$
\ms

Итак, для данного пакета $\calN$ мы получили набор множителей Лагранжа,
которые порождают функцию $\p(t),$ так что выполнены условия (\ref{stat1a}),
(\ref{stat2a}), (\ref{adj}), (\ref{tr1}), (\ref{usmax}).\, Эти множители
Лагранжа, вообще говоря, зависят от пакета.  Условия (\ref{stat1a}),
(\ref{stat2a}), (\ref{adj}), (\ref{tr1}) одни и те же для всех пакетов,
а условие максимума (\ref{usmax}) непосредственно связано с каждым отдельным
пакетом.  Наша цель теперь --- перейти к условию максимума для всех
$\theta \in (t_0, t_1)$ и $v \in U$ с набором множителей, не зависящим
от $\theta$ и $v.$
\ms

Предварительно заметим, что если $\alpha_0+|\alpha|+|\beta|=0,$ то $l=0,$
поэтому $\p(t_1)=0,$ и так как $\p$ удовлетворяет линейному однородному
уравнению, $\p(t) \equiv 0,$ и тогда в силу (\ref{lei}) получаем $\g =0,$
что противоречит условию нетривиальности в (\ref{stat2a}).\, Следовательно,
все определяется множителями $(\a_0, \a, \b),$ а условие нетривиальности
можно записать как $\alpha_0+|\alpha|+|\beta|>0.$ Заменим его условием
нормировки $\alpha_0+|\alpha|+|\beta|=1.$ \bs

%%-----------------------
{\bf Организация условий оптимальности, переход к универсальному ПМ.}\q
Обозначим через $\Lambda({\cal N})$  множество всевозможных наборов множителей
$(\alpha_0,\alpha,\beta),$ удовлетворяющих условиям
\begin{eqnarray} && \label{stat1aa}
\alpha_0\ge 0, \quad \alpha\ge 0,\quad \alpha_0+|\alpha|+|\beta|=1,
\quad \alpha F(\hat x_0,\hat x_1)=0, \\[4pt]
&& \label{stat2aa}
\dot \psi =\, -\psi f_x(\hat x(t),\hat u(t)), \qq
\p(t_0) = l_{x_0}\,, \q \p(t_1) = -l_{x_1}\,,
\end{eqnarray}
и конечнозначному условию максимума (\ref{usmax}). (Такой набор может быть
не один).\, Нетрудно видеть, что $\Lambda({\cal N})$ есть непустой компакт
в пространстве $\R^{1+d(F)+ d(K)}.$ \ssk

Таким образом, перебирая всевозможные пакеты $\calN,$ мы для каждого из них
получаем свой непустой компакт $\Lambda({\cal N}).\,$ Покажем, что семейство
всех этих компактов $\{ \Lambda({\cal N})\}_{\cal N}$ образует {\it
центрированную систему}. Для этого введем отношение порядка в множестве всех
пакетов. Будем говорить, что  ${\cal N}_1\subset{\cal N}_2\,,$ если каждая
пара $(\theta_i,v_i)$ из ${\cal N}_1$  входит (возможно с другим индексом $i$)
и в ${\cal N}_2\,.$
Ясно что для любых двух пакетов ${\cal N}_1$  и ${\cal N}_2$ найдется третий,
содержащий каждый из них, например их объединение. Далее ясно, что при расширении
пакета $\calN$ множество $\L(\calN)$ сужается, т.е. включение ${\cal N}_1\subset{\cal N}_2$
влечет обратное включение $\Lambda({\cal N}_1)\supset \Lambda({\cal N}_2).$
Пусть теперь имеется конечный набор компактов $\Lambda({\cal N}_1),\ldots,\Lambda({\cal N}_r)$.
Возьмем любой пакет ${\cal N},$ содержащий в себе все пакеты ${\cal N}_1,\ldots,{\cal N}_r\,.$
Тогда непустой компакт  $\Lambda({\cal N})$ содержится в каждом из компактов
$\Lambda({\cal N}_1),\ldots,\Lambda({\cal N}_r)$ и, следовательно, содержится
в их пересечении. Отсюда вытекает центрированность системы
$\{ \Lambda({\cal N})\}_{\cal N}\;$ и значит, непустота ее пересечения
$$
\L_*\; =\; \bigcap\limits_{{\cal N}}\,\Lambda({\cal N}).
$$

Возьмем произвольный набор множителей $(\alpha_0,\alpha,\beta)\in \L_*\,,$
и пусть $\psi(t)$ есть сопряженная функция, соответствующая этому набору.
Для этого набора по определению выполнены условия (\ref{stat1aa})--(\ref{stat2aa}),
а условие (\ref{usmax}) максимума функции Понтрягина по управлению выполнено
теперь {\it для любой пары} $(\theta,v) \in (t_0,t_1) \times U,$ ибо любая
такая пара входит в какой-нибудь пакет ${\cal N}.$ Из кусочной непрерывности
функции $\hu(t)$ следует, что тогда оно выполнено и для концов отрезка. \ssk

Таким образом, условия $(i')-(v')$ теоремы  \ref{th2} доказаны.
\bs

%%-----------------------------------
{\bf Постоянство функции Понтрягина для оптимального процесса.}\\
Нам осталось доказать, что вдоль оптимального процесса
\begin{equation}\label{const-h}
H(\psi(t),\hat x(t),\hat u(t))=\, \const \quad\mbox{ п.в.\, на } [t_0,t_1],
\end{equation}
т.е. выполнено условие $(vi').$ Тогда условие $(vii')$ будет следовать
из условия максимума (\ref{usmax}). Покажем, что условие (\ref{const-h})
вытекает из уже доказанных условий ПМ.\,\, Это можно сделать двумя способами:\,
1) показать это непосредственно,\,\, и\, 2) получить (\ref{const-h}) автоматически
из уже доказанных условий ПМ после перехода к задаче на нефиксированном
отрезке времени.\,

Второй способ, на наш взгляд, более простой и естественный (см. следующий
параграф). Достоинство первого способа в том, что он устанавливает
непосредственную зависимость условия (\ref{const-h}) от остальных условий
ПМ.\, Он рассматривался в книгах \cite{Pont, Bolt}.\, Мы здесь приводим
модификацию рассуждений из \cite{Bolt} (в книге \cite{Pont} проводятся
более сложные рассмотрения).
\ssk

Рассмотрим функцию $h(t,v) = H(\psi(t),\hat x(t),v) = \p(t)\,f(\hat x(t),v).$
Нам надо показать, что $h(t,\hu(t))$ постоянна на $[t_0,t_1].$
Будем рассуждать так.
\ssk

Во-первых, согласно условию максимума,
\beq{h1}
h(t,v) \lee h(t,\hu(t)) \qq \all t\in [t_0,t_1], \q \all v\in U.
\eeq

Далее, сама функция $h(t,v)$ и ее частная производная $h_t(t,v)$ непрерывны
на $[t_0,t_1]\times U,$ и в любой точке непрерывности управления $\hu(t)$
\beq{h2}
h_t(t,\hu(t))\;=\; 0.
\eeq
Действительно,
$$
h_t(t,\hu(t)) =\; \left.\Big(\dot\p(t)\,f(\hat x(t),v)\, +\,
\p(t)\,f_x(\hat x(t),v)\,\dot{\hx}(t)\Big)\right|_{v= \hu(t)} \; =
$$ $$
=\; \left.\Big(-\p(t)\,f_x(\hat x(t),\hu(t))\,f(\hat x(t),v)\,+\,
\p(t)\,f_x(\hat x(t),v)\,f(\hat x(t),\hu(t))\Big)\right|_{v=\hu(t)}\; =0.
$$

Покажем теперь, что в любой точке непрерывности $\hu(t)$ также и полная
производная $\frac d{dt}\, h(t,\hu(t))=0.$ Возьмем любую такую точку $t,$
и пусть точки $\t_2 > \t_1$ сходятся к $t.$ Положим $\Delta \t=\tau_2-\tau_1\,.$
Имеем
$$
h(\t_2, \hu(\t_2)) - h(\t_1,\hu(\t_1)) =\;\;
h(\t_2, \hu(\t_2))- h(\t_1, \hu(\t_2))  \;+ \phantom{pppppppp}
$$ \beq{ht12}
\phantom{ppppprrrrpyyy} +\; h(\t_1, \hu(\t_2)) - h(\t_1, \hu(\t_1)).
\eeq
В правой части этого равенства модуль первой разности оценивается по теореме
о конечном приращении:
$$
\Big|h(\t_2, \hu(\t_2))- h(\t_1, \hu(\t_2))\Big| \lee
\max_{\t\in [\t_1, \t_2]} \Big| h_t(\t,\hu(\t_2)) \Big|\cdot|\D\t| \; =\; o(\D\t),
$$
ибо стоящий здесь максимум стремится к $|\, h_t(t,\hu(t))| =0.\,$ Вторая
разность в (\ref{ht12}) неположительна в силу условия максимума (\ref{h1}).
Таким образом, $h(\t_2, \hu(\t_2)) - h(\t_1,\hu(\t_1))\le o(\D\t).$ \ssk

С другой стороны, эту же разность можно представить и так:
$$
h(\t_2, \hu(\t_2)) - h(\t_1,\hu(\t_1)) =\;\;
h(\t_2, \hu(\t_2))- h(\t_2, \hu(\t_1))  \;+ \phantom{pppppppp}
$$ \beq{ht21}
\phantom{ppppprrrrpyyy} +\; h(\t_2, \hu(\t_1)) - h(\t_1, \hu(\t_1)).
\eeq
откуда по тем же причинам вытекает, что
$\; h(\t_2, \hu(\t_2)) - h(\t_1,\hu(\t_1))\ge o(\D\t).$ \\
Следовательно,
$$
h(\t_2, \hu(\t_2)) - h(\t_1,\hu(\t_1))\; =\; o(\D\t).
$$
Отсюда вытекает, что полная производная $\dot h(t,\hu(t))$ существует и
равна нулю во всех точках непрерывности управления $\hat u(t),$ т.е.
$h(t,\hu(t)) = \const$ на каждом интервале непрерывности $\hat u(t).$
Покажем, что для всех интервалов эта константа одна и та же.
\ssk

Пусть $\tau\in (t_0,t_1)$ --- точка разрыва управления $\hat u(t)$
(согласно нашему предположению, их не более чем конечное число). Покажем,
что $h(t,\hu(t))$ не имеет разрыва в этой точке.\, При $t>\tau$ в силу
условия (\ref{h1}) имеем
$$
h(t,\hu(t)) \gee h(t,\hu(\t-0)).
$$
Переходя к пределу по $t\to\tau+0\,,$ получаем
$$
h(\t,\hu(\t+0)) \gee h(\t,\hu(\t-0)).
$$
Аналогично при $t<\tau$ в силу (\ref{h1}) имеем
$$
h(t,\hu(t)) \gee h(t,\hu(\t+0)).
$$
Переходя к пределу по $t\to\tau-0\,,$ получаем
$$
h(\t,\hu(\t-0)) \gee h(\t,\hu(\t+0)).
$$
Следовательно, $h(\t,\hu(\t-0)) = h(\t,\hu(\t+0)).$
Отсюда и из предыдущего вытекает, что $h(t,\hu(t))= \const$ на всем отрезке
$[t_0, t_1].$  Тем самым условие (\ref{const-h}) доказано, а вместе
с ним полностью доказана и теорема~\ref{th2}.
\bs

%%------------------------------------------------
\section{Доказательство принципа максимума для\\ задачи на нефиксированном отрезке времени}

Вернемся к задаче (\ref{4})--(\ref{1}) на нефиксированном отрезке времени.
Мы назвали ее задачей A.\, Пусть
$ w^0=(x^0(t), u^0(t))\mid  t\in [\hat t_0, \hat t_1])$  ---
допустимый процесс в этой задаче, для которого имеет место сильный минимум.
(Здесь нам будет удобно снабжать $x,\,u$ индексом 0.)\,
Как и ранее, предположим, что управление  $ u^0(t)$ кусочно непрерывно.\,

Сведем задачу А\, к некоторой вспомогательной задаче на {\em фиксированном}\,
отрезке времени, с тем, чтобы затем воспользоваться теоремой \ref{th2}.
Для этого введем "новое время" $\,\tau,$ меняющееся на отрезке
$[\tau_0, \tau_1]:=[\hat t_0,\hat t_1],$ а исходное время $t$ будем рассматривать
как еще одну фазовую переменную $t= t(\t),$ подчиненную уравнению
\beq{dttau}
\frac {dt}{d\t}\; =\; v(\t),
\eeq
где $v(\t)>0$ есть еще одно управление.
%Последнее неравенство вместе с исходным множеством $\calR,$ задает открытое
%множество $\calR' = \calR\times \R_+ \subset \R^{2+n+r},$ on which the right
%hand side of the new control system is defined and sufficiently smooth.
Исходные переменные $x,\, u$ в новом времени есть $\tl x(\t) = x(t(\t)),$
$\tl u(\t)= u(t(\t));\,$ они подчинены уравнению
\beq{eq1tau}
\frac {d\tl x}{d\t}\; =\; v(\t)\,f(t(\t), \tl x(\t), \tl u(\t)).
\eeq
Концевой блок задачи выглядит теперь так:
\beq{eq2tau}
K(t(\t_0),\tl x(\t_0), t(\t_1), \tl x(\t_1)) =0, %\qq j=1,\ld, k,
\eeq
\beq{eq3tau}
F(t(\t_0),\tl x(\t_0), t(\t_1), \tl x(\t_1)) \le 0, %\qq i=1, \ld, \nu,
\eeq
\beq{eq4tau}
J = F_0(t(\t_0),\tl x(\t_0), t(\t_1), \tl x(\t_1)) \to \min.
\eeq \\[2pt]
Мы видим, что он зависит только от концевых значений фазовых переменных $t$
и $\tl x.$ Задачу (\ref{dttau})--(\ref{eq4tau}) назовем задачей $\wt A.$
% а переменные $x,\,u$ в ней будем снабжать тильдой.
В ней $t(\tau),\;\wt x(\tau)\;-$ фазовые переменные, а
$v(\tau),\; \wt u(\tau)\;-$ управления. Условие $v>0$ следует рассматривать
как ограничение на новое управление. \ssk

Процессу $w^0$  в  задаче A\, мы сопоставим процесс
$$
\wt w^0\,=\;( v^0(\tau),\, t^0(\tau),\, \wt x^0(\tau),\, \wt u^0(\tau))\;
\mid \; \tau\in [\tau_0, \tau_1]),
$$
в  задаче $\wt A$, где
\begin{equation}\label{3nn}
v^0(\tau)\equiv 1, \q t^0(\tau)\equiv\tau, \quad
\wt x^0(\tau)\equiv x^0(\tau), \quad
\wt u^0(\tau)\equiv u^0(\tau).
\end{equation}
Ясно, что этот процесс допустим в задаче $\wt A.$ Более того, из того, что
процесс $w^0$ доставляет сильный минимум в задаче A, вытекает, что процесс
$\wt{w}^0$ доставляет сильный минимум в задаче $\wt A.$ Это легко устанавливается
от противного. По теореме~\ref{th2}\, для $\wt {w}^0$ выполнен принцип
максимума. Выпишем его условия для процесса $\wt {w}^0$ в задаче $\wt A,$
опуская теперь тильду, и учитывая, что $dt^0/ d\t \equiv 1$ и
$t^0(\t) \equiv \t,$ так что производную по $\t$ по-прежнему можно
обозначать точкой.
\ms

Концевая функция Лагранжа в задаче $\wt A$ имеет тот же вид, что и в задаче A:
$$
l=\; \alpha_0 F_0+\alpha F+\beta K.
$$
Функция Понтрягина в задаче $\wt A$ имеет вид
$$
\wt H(\p_x,\p_t,\,t,x,u,v) =\; (\p_x\,f(t,x,u) + \p_t)\,v \; =
\big(H(\p_x,t,x,u) + \p_t \big)\,v.
$$

Условия принципа максимума для процесса $w^0$  в задаче $\wt A$ таковы:
\ssk

\noi
существует число $\alpha_0\,,$ векторы $\alpha\in\R^{d(F)},$ $\beta\in \R^{d(K)},$
липшицевы функции $\psi_t(\t)$ и $\psi_x(\t)$ размерности 1\, и $n$
соответственно, такие что выполнены условия:
\ssk

\noi
неотрицательности и нетривиальности
\begin{equation}\label{mp1}
\alpha_0\ge 0,\quad \alpha\ge 0; \quad \; \alpha_0+|\alpha|+|\beta|>0,
\eeq
дополняющей нежесткости
\beq{dopnez}
\alpha F(t^0(\tau_0), x^0(\tau_0), t^0(\tau_1), x^0(\tau_1))=0,
\end{equation}
сопряженные уравнения
\beq{dpsixtau}
-\dot\psi_x(\t) = \wt H_x \;=\; \p_x(\t)\,f_x(t^0(\t),x^0(\t),u^0(\t)),
\eeq
\beq{dpsittau}
-\dot\psi_t(\t) = \wt H_t \;=\; \p_x(\t)\,f_t (t^0(\t),x^0(\t),u^0(\t)),
\eeq
трансверсальности
\beq{transxtau}
\p_x(\t_0) =\; l_{x(\t_0)}\,, \qq  \p_x(\t_1) =\; -l_{x(\t_1)}\;,
\eeq
\beq{transttau}
\p_t(\t_0) =\; l_{t(\t_0)}\,, \qq  \p_t(\t_1) =\; -l_{t(\t_1)}\;,
\eeq
и условие максимальности по управлениям $u, v$:
\begin{equation}\label{mp10}
v\Big(\psi_x(\tau)\,f(t^0(\tau), x^0(\tau),u)+ \psi_t(\tau)\Big)\le
v^0(\tau)\Big(\psi_x(\tau)\,f(t^0(\tau), x^0(\tau), u^0(\tau))+\psi_t(\tau)\Big)
\eeq
$$
\mbox{для всех}\q u\in U, \q v>0,\;\q \tau\in[\tau_0,\tau_1].
$$ \ssk

\noi
При $u = u^0(\t)$ из максимальности  $\wt H$ по $v>0$ в точке $v^0 =1$
следует, что
\begin{equation}\label{mp13}
\psi_x(\tau)\,f(t^0(\tau), x^0(\tau), u^0(\tau)) +\psi_t(\tau)\,= 0
\quad\; \mbox{ на }\; [\tau_0,\tau_1],
\end{equation}
и тогда (\ref{mp10}) означает, что
\begin{equation}\label{mp14}
\psi_x(\tau)\,f(t^0(\tau), x^0(\tau), u^0(\tau))+\psi_t(\tau)\, \le 0\quad
\mbox{для всех}\quad u\in U, \quad \tau\in[\tau_0,\tau_1].
\end{equation}
Учитывая, что $t^0(\tau) =\t,$ мы видим, что получены все условия
$(i)-(vii)$ принципа максимума для процесса $w^0$ в задаче A.\,
Теорема~\ref{th1}\, полностью доказана. \ctd \ms

Отметим, что для автономной задачи из (\ref{dpsittau}) вытекает, что
$\p_t(\t) = \const,$ и тогда (\ref{mp13}) означает, что
$H(\p_x(t), x^0(t), u^0(t)) = -\p_t = \const,$ т.е., как было сказано выше,
мы автоматически получаем условие (\ref{const-h}).  \bs

%%-------------------------------
\section{Приложение}

Для доказательства леммы \ref{lem2}\, нам понадобится следующая
\ble{lem7}
Пусть на отрезке $[a,b]$ задана непрерывная функция $f,$ такая что для
всех $x \in [a,b)$ существует $f'_{\mbox{\rm \small пр}}(x) >0.\;$ Тогда
$f(b) > f(a).$ Отсюда также следует, что $f$ строго возрастает на $[a,b].$
\ele

\Proof Считаем $f(a)=0.$ Так как $f'_{\mbox{\rm \small пр}}(a) >0,$ то
$f>0$ на интервале $(a,\,a+\d)$ при некотором $\d>0.$ Возьмем любую точку
$a_1$ из этого интервала, и пусть $a_2$ есть максимальная из тех точек $x,$
для которых $f \ge f(a_1)$ на отрезке $[a_1,\,x].$ Если $a_2 <b,$ то поскольку
$f'_{\mbox{\rm \small пр}}(a_2) >0,$ %то
в некоторой правой полуокрестности
точки $a_2$ будет $f > f(a_2) \ge f(a_1),$ что противоречит максимальности $a_2\,.$
Следовательно, $a_2 =b,$ и тем самым лемма доказана. \ctd \bs

С помощью этой леммы докажем аналог теоремы Лагранжа о конечном приращении
для правой производной.

\ble{lem8}
Пусть функция $f$ непрерывна на отрезке $[a,b],$ и для всех $x \in [a,b)$
имеет $f'_{\mbox{\rm \small пр}}(x),$ непрерывную по $x.$ Тогда
$$
f(b) - f(a)\, =\, f'_{\mbox{\rm \small пр}}(\tl x)\,(b-a) \qq
\mbox{при некотором}\q \tl x\in [a,b).
$$
\ele

\Proof Вычитая из $f$ линейную функцию, приходим к случаю $f(a) = f(b) =0$
(т.е. надо доказать аналог теоремы Ролля). Допустим, что утверждение леммы
неверно, например $f'_{\mbox{\rm \small пр}}(x) >0$ для всех $x \in [a,b).$
Тогда по лемме \ref{lem7}\, имеем $f(b) > f(a),$ что противоречит условию.
\ctd \bs

{\bf Доказательство леммы \ref{lem4}}.\, Возьмем любую точку $x >0.$
Без нарушения общности считаем, что $\f'_{\mbox{\rm \small пр}}(x) =0.$
Нам надо показать, что и левая производная
$\f'_{\mbox{\rm \small л}}(x) =0,$ т.е. что
$|\,\f(x) - \f(x-\d)| = o(\d)$ при $\d \to 0+.$ По лемме \ref{lem8}}\,
$$
|\,\f(x) - \f(x-\d)|\; =\; \f'_{\mbox{\rm \small пр}}(\tl x)\,\d \qq
\mbox{при некотором} \q \tl x\in [x,\,x-\d).
$$
В силу непрерывности правой производной
$\f'_{\mbox{\rm \small пр}}(\tl x) \to \f'_{\mbox{\rm \small пр}}(x) =0,$
откуда следует требуемое. \ctd \bs

{\bf Доказательство леммы \ref{lem5}}.\, Не нарушая общности считаем,
что $P'(0,0)=0,$ т.е. отображение $P$ липшицево в $r-$окрестности точки
$(0,0) \in \R^n \times K$ с константой $\mu(r)\to 0$ при $r \to 0.$
Зафиксируем любой вектор $h \in \inter K,$ и для каждого $y \in \R^s$
положим $\l(y) = \min\,\{\l \ge 0\;|\; y +\l h \in K\}.$  Построим
отображение $\f:\, \R^s \to K$ по формуле $\f(y) = y + \l(y)\,h.$
(Это есть проекция $\R^s$ на $K$ вдоль вектора $h.)$ Для всех $y\in K$
очевидно $\f(y) =y.$ Ясно, что $\f$ липшицево с некоторой константой $L.$
Тогда отображение $\wt P(x,y) = P(x,\,\f(y))$ определено на всем
$\R^n \times \R^s$ и липшицево в $r-$окрестности точки $(0,0)$ с константой
$L\,\mu(r)\to 0.$ Следовательно, оно имеет строгую производную в нуле,
равную нулю. \ctd
\ms

Отметим, что эта лемма справедлива и в случае более общего отображения
$P: X \times K \to Z,\,$ где $X,\, Y,\, Z\;-$ произвольные нормированные пространства,
а $K \subset Y$ есть замкнутый выпуклый конус с непустой внутренностью.
Доказательство остается без изменений.\,
%В случае $K = \R^s_+$ можно положить $h = (1, \ld, 1),$
%и тогда $\l(y) = -\min\,\{y_1, \ld, y_s\}.$
\ms

Для доказательства леммы \ref{lem5a}, как и в \cite{Kor}, установим следующее,
%хорошо известное (см. напр. \cite{Lich})\, и по сути одномерное свойство.
по сути одномерное (и конечно хорошо известное, см. напр. \cite{Lich}), свойство.

\ble{lem9}
Пусть отображение $f: \R^n \times \R_+ \times \R^k_+ \to \R^m$ имеет производную
во всех внутренних точках области определения, которая непрерывна вплоть до
границы этой области. Продолжим $f$ на $\R^n \times \R \times \R^k_+\,,$
положив для $y<0$
$$
\wt f(x,y,z)\, =\; -f(x, -y,z)\,+\, 2 f(x,0,z).
$$
Тогда полученное отображение $\wt f$ обладает теми же свойствами гладкости
уже на своей области определения.
\ele

Доказательство представляет собой элементарную проверку.
\ssk

{\bf Доказательство леммы \ref{lem5a}}\, состоит
%Лемма \ref{lem5a}\, вытекает
из последовательного применения
леммы~\ref{lem9}\, ко всем компонентам вектора $\e.$  \\

Авторы выражают благодарность\, В.А. Дыхте\, за ценные замечания\, и \\
Г.Г. Магарил--Ильяеву за полезные обсуждения.
\bs

\small
%\np
%%-------------------------------------------------------------------

%%--------------------------------------------------------------
\vskip 8mm

{\bf Дмитрук Андрей Венедиктович,}\, ЦЭМИ РАН,\, Москва 117418, Нахимовский
проспект, 47;\q МГУ им. Ломоносова, мехмат\,  %Тел. +7-909-940-87-13,\\
({\tt dmitruk@member.ams.org}). \ms

{\bf Осмоловский Николай Павлович,}\, Московский госу\-дарст\-венный
строительный университет,\, Москва, Ярославское шоссе, 26;\q
Университет технических и гуманитарных наук г. Радома, Польша,\,
%ул. Мальчевского, 20-а \,  \quad Тел. +7-903-249-63-98, \q
({\tt osmolovski@uph.edu.pl}).

\normalsize